\documentclass[a4paper,12pt,reqno]{amsart}
\usepackage[T1]{fontenc}
\usepackage[utf8]{inputenc}
\usepackage[english]{babel}
\usepackage{amsthm,amsfonts,amsmath,amssymb}
\usepackage[left=3.5cm,right=3.5cm,top=3.5cm,bottom=3.5cm,heightrounded]{geometry}
\usepackage[dvipsnames]{xcolor}
\usepackage{mathtools}
\usepackage{enumitem}

\numberwithin{equation}{section}
\allowdisplaybreaks

\newcommand{\R}{\mathbb{R}}

\DeclareMathOperator{\dist}{dist}
\DeclareMathOperator{\Per}{Per}
\DeclareMathOperator{\Div}{div}
\DeclareMathOperator{\co}{co}
\DeclareMathOperator{\ri}{ri}

\usepackage[colorlinks=true,urlcolor=teal,citecolor=magenta,linkcolor=teal,linktocpage,pdfpagelabels,bookmarksnumbered,bookmarksopen,breaklinks=true]{hyperref}
\usepackage[capitalize,nameinlink,noabbrev]{cleveref}

\theoremstyle{plain}
	\newtheorem{theorem}{Theorem}[section]
    \newtheorem{lemma}[theorem]{Lemma}
    \AddToHook{env/lemma/begin}{\crefalias{theorem}{lemma}}
    \crefname{lemma}{Lemma}{Lemmas}
	\newtheorem{proposition}[theorem]{Proposition}
    \AddToHook{env/proposition/begin}{\crefalias{theorem}{proposition}}
    \crefname{proposition}{Proposition}{Propositions}
	\newtheorem{corollary}[theorem]{Corollary}
    \AddToHook{env/corollary/begin}{\crefalias{theorem}{corollary}}
    \crefname{corollary}{Corollary}{Corollaries}
	
    \newtheorem*{mainthm}{Main Theorem}
    \newcounter{mainthm}
    \AddToHook{env/mainthm/begin}{\refstepcounter{mainthm}}
    \crefformat{mainthm}{#2Main Theorem#3}
    \Crefformat{mainthm}{#2Main Theorem#3}

     \newtheorem*{maincor}{Corollary}
    \newcounter{maincor}
    \AddToHook{env/maincor/begin}{\refstepcounter{maincor}}
    \crefformat{maincor}{#2Corollary#3}
    \Crefformat{maincor}{#2Corollary#3}
    
    \theoremstyle{definition}
	
    \AddToHook{env/definition/begin}{\crefalias{theorem}{definition}}
    \crefname{definition}{Definition}{Definitions}
	\newtheorem{remark}[theorem]{Remark}
    \AddToHook{env/remark/begin}{\crefalias{theorem}{remark}}
    \crefname{remark}{Remark}{Remarks}
	\newtheorem{example}[theorem]{Example}
    \AddToHook{env/example/begin}{\crefalias{theorem}{example}}
    \crefname{example}{Example}{Examples}

\title[Admissible geometries for Talenti's P\'olya--Szeg\H{o} principle]{Admissible geometries \\
for Talenti's P\'olya--Szeg\H{o} principle}

\author[G.~Bartoli]{Giulio Bartoli}
\address[G.~Bartoli]{
Institut f\"ur Mathematik, Goethe-Universit\"at Frankfurt, \linebreak Robert-Mayer-Str.~10, 60325 Frankfurt am Main, Germany
}
\email{bartoli@mathematik.uni-frankfurt.de}

\author[G.~Saracco]{Giorgio Saracco}
\address[G.~Saracco]{
Dipartimento di Matematica e Informatica, Universit\`a di Ferrara, via Machiavelli 30, 44121 Ferrara (FE), Italy
}
\email{giorgio.saracco@unife.it}



\begin{document}

\begin{abstract}
Talenti symmetrization rearranges functions on an open connected set \(X\subset\R^N\), endowed with a measure \(\mu\) with continuous positive density, along the superlevel sets of a given function \(f\).
Talenti proved a P\'olya--Szeg\H{o} principle for it, provided that these sets are \(\mu\)-isoperimetric, have connected boundaries, and are crossed by straight lines of steepest ascent.
We show that, if \(X\) is convex, these requirements are very rigid: they force the superlevel sets of \(f\) to be either parallel half-spaces, or concentric balls, or complements of tubes around an affine subspace of codimension at least two, and, if \(\mu(X)<+\infty\), only parallel half-spaces occur.
This identifies the only geometries in which Talenti's P\'olya--Szeg\H{o} principle can be applied.
We also prove a family of weighted isoperimetric inequalities in convex cones times cubes, showing that complements of tubes around affine subspaces of every codimension at least two do occur.
\end{abstract}

\maketitle
\section{Introduction}
\label{sec:intro}
Let \(X\subset\R^N\), with \(N\ge2\), be an open connected set, and let \(\mu=e^W\mathcal{L}^N\), with \(W\) continuous on \(X\).
In~\cite{Tal97}, Talenti introduced a symmetrization which rearranges a function \(u\) along the superlevel sets of a fixed function \(f\): the symmetrized function \(u^\sharp\) is a nonincreasing function of \(f\), \(\mu\)-equimeasurable with \(u\).
In analogy with the classical P\'olya--Szeg\H{o} inequality~\cite{PS51book}, which states that the symmetric decreasing rearrangement does not increase the \(L^p\) norm of the gradient, Talenti proved that his symmetrization does not increase the \(L^p(\mu)\) norm of the gradient either, for bounded Lipschitz functions \(u\) with \(\mu(\{u\neq0\})<+\infty\), see~\cite[Main Thm.]{Tal97}.
In~\cite{BS26}, assuming \(X\) convex, we extended this result to Sobolev functions vanishing on the relative boundary of an open set \(\Omega\subset X\), thus obtaining Faber--Krahn and Saint-Venant inequalities in a variety of weighted settings, see~\cite[Sect.~6]{BS26} and the references therein.

Talenti's assumptions, recalled in \cref{thm:talenti} below, are that the superlevel sets of \(f\) are isoperimetric in \((X,\mu)\), relative to \(X\), that the level sets of \(f\) are connected, and that the lines of steepest ascent of \(f\) are straight, see~\ref{hp:T1}--\ref{hp:T4} in \cref{sec:setting}.
In addition, \(f\) is required to be smooth and satisfy some normalization and nondegeneracy conditions, see~\ref{hp:N1}--\ref{hp:N3} in \cref{sec:setting}.
Neither the measure nor the function \(f\) is prescribed, and this generality is what makes the principle attractive: Schwarz symmetrization~\cite{Sch90book}, with balls and the Lebesgue measure, and Ehrhard symmetrization~\cite{Ehr83}, with half-spaces and the Gaussian measure, are particular instances.

It is then natural to ask which geometries are admissible, that is, which families of sets can arise as the superlevel sets of \(f\).
Since these sets are isoperimetric in \((X,\mu)\), relative to \(X\), the question becomes which families of isoperimetric sets can occur in weighted settings, as the density varies, while the other two assumptions, connectedness of the level sets and straightness of the lines of steepest ascent, still hold.
In all the examples we are aware of, including those of~\cite[Sect.~6]{BS26}, these sets are either parallel half-spaces, or concentric balls, or complements of balls.
The purpose of this paper is to show that, when \(X\) is convex, this is not a coincidence: Talenti's assumptions single out these three families, up to allowing, in the third one, tubes around affine subspaces in place of balls.
Our main result reads as follows.

\begin{mainthm}
\label{thm:main}
Let \(X\subset\R^N\), with \(N\ge2\), be an open convex set, let \(\mu=e^W\mathcal{L}^N\), with \(W\) continuous on \(X\), and let \(f\colon X\to[0,1]\) be a smooth function satisfying~\ref{hp:N1}--\ref{hp:N3} and~\ref{hp:T1}--\ref{hp:T4}.
Then, the level sets of \(f\) are parallel hypersurfaces, and, relative to \(X\), its superlevel sets are either
\begin{enumerate}[label=(\alph*)]
\item parallel half-spaces, or
\item concentric balls, centered at a point of \(\overline{X}\), or
\item complements of closed tubes around an affine subspace \(L\), with \(\dim L\le N-2\) and \(L\cap X=\emptyset\), which is the affine hull of \(L\cap\overline{X}\).
\end{enumerate}
Moreover, if \(\mu(X)<+\infty\), only~(a) may occur, whereas alternative~(c) requires \(\mu(B_\rho(x)\cap X)=+\infty\) for all \(x\in L\cap\overline{X}\) and \(\rho>0\).
\end{mainthm}

We single out the case of finite measure, which shows that, on convex sets of finite measure, Talenti symmetrization is necessarily a symmetrization with respect to parallel half-spaces, as Ehrhard's, the measure being the only remaining freedom.

\begin{maincor}
Under the assumptions of the \cref{thm:main}, if \(\mu(X)<+\infty\), then the superlevel sets of \(f\) are parallel half-spaces intersected with \(X\).
\end{maincor}

In other words, the \cref{thm:main} is a rigidity result: on convex ambient sets, every Faber--Krahn or Saint-Venant type inequality obtained through Talenti symmetrization, such as those of~\cite{BS26}, has parallel half-spaces, concentric balls, or complements of tubes as extremal sets.
Hence, the search for new applications essentially reduces to finding measures for which these sets are isoperimetric and have connected boundaries relative to \(X\).

This rigidity is of a different nature than the characterization of equality cases in P\'olya--Szeg\H{o} inequalities, see, \emph{e.g.},~\cite{BZ88}, and than stability results for symmetrization inequalities, see, \emph{e.g.},~\cite{AL26} for Talenti's comparison principle: rather than the extremal functions, it constrains the frameworks in which the inequality can be established.
It is closer in spirit to the classification of isoparametric and transnormal functions~\cite{Bol73,Seg38,Wan87}, from which it differs in two respects: the presence of a boundary and the isoperimetric property, which restricts the level sets of \(f\), a priori only parallel hypersurfaces, to be isoparametric, that is, hyperplanes, spheres, and cylinders by Segre's classification~\cite{Seg38}.

Two further features of the \cref{thm:main} are worth stressing.
First, the classification is not invariant under complementation: complements of tubes around affine subspaces of every dimension up to \(N-2\) occur, whereas tubes only occur around points, that is, as balls.
Tubes around affine subspaces of positive dimension are excluded by the isoperimetric property of the superlevel sets of small measure, even though they may be isoperimetric for large volumes~\cite{FMP11,RV15b}, see \cref{rem:focal}.
Second, Talenti's assumptions constrain the measure: besides being infinite in alternatives~(b) and~(c), by the \cref{thm:main}, in alternative~(c) the density cannot be integrable near any point of \(L\cap\overline{X}\), so that it necessarily blows up on the whole set \(L\cap\overline{X}\subset\partial X\).
In a complementary direction, fixing the geometry constrains the measure even further: for instance, if half-spaces are isoperimetric in a cylinder, then the density has a product structure, as shown by Brock, Chiacchio, and Mercaldo~\cite{BCM16}.

All the alternatives of the \cref{thm:main} do occur.
Concentric balls and parallel half-spaces appear in~\cite[Sects.~6.1 and~6.2]{BS26}, see also~\cite{BCM12,BCM14,CROS16} for balls centered at the vertex of a convex cone and~\cite{BCM16} for half-spaces with Gauss-like product measures, while the complements of the disks centered at the vertex of a planar convex sector, endowed with a density blowing up at the vertex, provide an instance of alternative~(c) with \(\dim L=0\), see~\cite[Sect.~6.3]{BS26} and~\cite{CJQW08,DHHT12}.
In \cref{sec:examples-c} we prove, via calibration, a family of weighted isoperimetric inequalities in convex cones times cubes, which generalizes this planar example and may be of independent interest, thus showing that alternative~(c) occurs with \(\dim L=\ell\) for every \(\ell\in\{0,\dots,N-2\}\).
For \(\ell=N-1\), the same construction yields parallel half-spaces in a setting with \(\mu(X)=+\infty\), so that alternative~(a), unlike~(b) and~(c), occurs both for finite and for infinite measures, see \cref{rem:hyperplane}.
The same calibration argument, applied in the exterior of a convex body, also shows that convexity cannot be dropped in the \cref{thm:main}, see \cref{rem:nonconvex}.

\subsection*{Outline of the paper and strategy of the proof}
In \cref{sec:setting}, we recall the setting, the normalization and nondegeneracy conditions~\ref{hp:N1}--\ref{hp:N3}, Talenti's assumptions~\ref{hp:T1}--\ref{hp:T4}, and the results of~\cite{BS26,Tal97} that we need.

In \cref{sec:structure}, we describe the geometry of \(f\): by~\cite[Thm.~4.1]{BS26}, \(\|\nabla f\|=g(f)\) for a continuous profile \(g\), and, the level sets being connected, the straightness of the lines of steepest ascent amounts to \(f\) being, away from its extremal level sets, \emph{transnormal} in the sense of~\cite{Bol73,Wan87}.
Unlike in the classical setting, where the ambient space has no boundary, the lines of steepest ascent may leave \(X\), and the focal sets of \(f\) typically lie on \(\partial X\), so that the known structure results do not apply, see \cref{rem:transnormal}.
We show that, if \(f\) satisfies a \emph{sweeping condition}, that is, if it takes all its values along each of its lines of steepest ascent, then it is a function of the signed distance from any of its level sets, and either its level sets are parallel hyperplanes, or \(f\) is a function of the distance from a focal set, see \cref{thm:dist-rep,cor:parallel}.

In \cref{sec:isoperimetry}, we combine this structure with the isoperimetric property.
The isoperimetric arguments rely on a single comparison principle, see \cref{lem:ode}: parametrized by arclength along the lines of steepest ascent, the measure \(V\) of the superlevel sets solves \(V'=-I(V)\), where \(I\) denotes the \(\mu\)-isoperimetric profile.
On the other hand, the measure \(w\) of any family of sets whose boundaries move at unit normal speed, such as parallel sets or concentric balls, satisfies \(|w'|\ge I(w)\), so that \(w\) can be compared with \(V\).
Applied to the inner and outer parallel sets of a superlevel set, it shows that the isoperimetric property forces the sweeping condition above, see \cref{lem:sweep}.
Applied to balls centered at points of \(\overline{X}\), it proves that each such ball has measure at least that of the tube of the same radius around a focal set.
Hence, if the superlevel sets of small measure are tubes around a focal set, they cannot contain disjoint balls centered at two distinct focal points; therefore, the focal set is a point, and the superlevel sets are concentric balls, see \cref{prop:focal-point}.
When \(\mu(X)<+\infty\), we prove the isoperimetric profile to be symmetric, which, together with the previous step, rules out focal sets, so that only parallel half-spaces occur, see \cref{prop:width}.
Applied to the complements of the metric neighborhoods of closed subsets of the focal set, the comparison principle shows that the density is not integrable near any of its points, see \cref{prop:non-integrable}.
In alternative~(c), the lines of steepest ascent are half-lines emanating from the focal set \(F_0\): this forces \(\overline{F_0}\) to be convex and, by the non-integrability of the density, to lie on \(\partial X\), while the continuity of the second fundamental form of the level sets of \(f\) forces \(\overline{F_0}\) to be the intersection of \(\overline{X}\) with an affine subspace, see \cref{prop:focal-c}.
These steps are combined in the proof of the \cref{thm:main}.

Finally, in \cref{sec:examples-c}, we prove the isoperimetric inequalities in cones briefly described above.
Through~\cite{BS26}, they yield Faber--Krahn, Saint-Venant, and \((p,q)\)-eigenvalue inequalities whose extremal sets are complements of tubes, see \cref{ex:wedge}.

\subsection*{Acknowledgments}
Both authors are members of INdAM--GNAMPA.
G.B.\ is supported by the Deutsche Forschungsgemeinschaft (DFG, German Research Foundation), Project No.~555837013.


\subsection*{Data availability statement}

This manuscript has no associated data.


\section{Setting and preliminaries}
\label{sec:setting}

\subsection{Assumptions}
Throughout the paper, \(X\subset\R^N\), with \(N\ge2\), is an open convex set, and \(\mu=e^W\mathcal{L}^N\), with \(W\) continuous on \(X\).
Convexity is the setting of~\cite{BS26}, where it is used to extend the P\'olya--Szeg\H{o} principle to Sobolev functions vanishing on the relative boundary of open subsets of \(X\), see~\cite[Thm.~5.1]{BS26}.
Here, it is used from \cref{thm:geo-uc} on, where segments joining points of \(X\) are required to lie in \(X\).
Talenti's theorem, recalled in \cref{thm:talenti}, only requires \(X\) to be open and connected, but convexity cannot be dropped in the \cref{thm:main}, see \cref{rem:nonconvex}.

As \(W\) is continuous, \(\mu\) and \(\mathcal{L}^N\) are mutually absolutely continuous, so that \(\mu(A)>0\) for every nonempty open set \(A\subset X\).
The \(\mu\)-perimeter relative to \(X\) of a set \(E\subset X\) of locally finite perimeter in \(X\) is
\begin{equation}
\label{eq:P-mu}
\Per_\mu(E;X)=\int_{\partial^*E\cap X}e^{W}\,\mathrm{d}\mathcal{H}^{N-1},
\end{equation}
where \(\partial^*E\) is the reduced boundary of \(E\); in particular, \(\partial X\) is not counted.
If a Borel set \(E\subset X\) does not have locally finite perimeter in \(X\), we let \(\Per_\mu(E;X)=+\infty\).
For the theory of sets of finite perimeter, we refer the interested reader to the comprehensive treatise~\cite{Mag12book}.

Throughout the paper, \(f:X\to[0,1]\) is a fixed smooth function; for \(t\in[0,1]\), we let \(\Phi_f(t) = \mu(\{f>t\})\).
We collect here all the assumptions on \(f\) and \(\mu\) that we shall use; each statement below recalls the ones it needs.

The first three are normalization and nondegeneracy conditions:
\begin{enumerate}[label=\textup{(N\arabic*)}, ref=\textup{(N\arabic*)}, leftmargin=*, widest=3]
\item
\label{hp:N1}
\(\inf f=0\) and \(\sup f=1\), with \(\mu\bigl(f^{-1}(0)\bigr)=\mu\bigl(f^{-1}(1)\bigr)=0\);
\item
\label{hp:N2}
\(\nabla f(x)\neq0\) for all \(x\in X\) with \(0<f(x)<1\);
\item
\label{hp:N3}
\(\Phi_f(t)<+\infty\) for all \(t\in(0,1)\), or \(\mu(X)<+\infty\).
\end{enumerate}
We let \(X_{\mathrm{reg}}=\{\,x\in X\,:\,0<f(x)<1\,\}\) and, for \(t\in(0,1)\), we denote by \(M_t=\{\,x\in X\,:\,f(x)=t\,\}\) the level sets of \(f\), which, under~\ref{hp:N2}, are smooth hypersurfaces.
Further, the unit normal
\begin{equation}
\label{eq:def-N}
\mathbf{N}(x)=\frac{\nabla f(x)}{\|\nabla f(x)\|},
\qquad x\in X_{\mathrm{reg}},
\end{equation}
is smooth on \(X_{\mathrm{reg}}\), and the \emph{lines of steepest ascent} of \(f\) are the maximal integral curves of \(\mathbf{N}\) in \(X_{\mathrm{reg}}\), that is, the maximal solutions \(s\mapsto x(s)\) of \(\mathrm{d}x/\mathrm{d}s=\mathbf{N}(x(s))\), where \(s\) is the arclength parameter.
The next four are the assumptions of Talenti's theorem:
\begin{enumerate}[label=\textup{(T\arabic*)}, ref=\textup{(T\arabic*)}, leftmargin=*, widest=4]
\item\label{hp:T1}
\emph{isoperimetric inequality:} there exists a continuous function \(I\ge0\) on \([0,\mu(X))\), positive on \((0,\mu(X))\), such that
\[
I(\mu(E))\le\Per_\mu(E;X),
\quad
\text{for all Borel sets \(E\subset X\) with \(\mu(E)<\mu(X)\);}
\]
\item\label{hp:T2}
\emph{saturation:} \(\Per_\mu(\{f>t\};X)=I(\Phi_f(t))\) for all \(t\in(0,1)\);
\item\label{hp:T3}
\emph{connectedness:} the level sets \(M_t\), \(t\in(0,1)\), are connected;
\item\label{hp:T4}
\emph{straightness:} the lines of steepest ascent of \(f\) are segments, that is, \(\mathrm{d}^2x/\mathrm{d}s^2=0\) along each of them.
\end{enumerate}
Finally, we shall use the following \emph{sweeping condition}, which is not among Talenti's assumptions, but follows from them, see \cref{lem:sweep}:
\begin{enumerate}[label=\textup{(SW)}, ref=\textup{(SW)}, leftmargin=*, widest=4]
\item\label{hp:SW}
along each line of steepest ascent of \(f\), the function \(f\) is onto \((0,1)\).
\end{enumerate}

\subsection{Known results}
The \((\mu,f)\)-Talenti symmetrization of a function \(u\) is defined by \(u^\sharp(x)=u^*(\Phi_f(f(x)))\), where \(u^*\) is the \(\mu\)-decreasing rearrangement of \(u\), see~\cite{Tal97}.
We shall only use the following properties of \(f\), see~\cite[Lem.~3.5 and~3.8]{BS26}.

\begin{lemma}
\label{lem:Phi-f}
Assume~\ref{hp:N1}--\ref{hp:N3}.
Then, \(\Phi_f\) is a continuous and strictly decreasing map of \([0,1]\) onto \([0,\mu(X)]\) (onto \([0,+\infty]\), in the case \(\mu(X)=+\infty\)).
Moreover, for every \(t\in(0,1)\), the level set \(M_t\) is a smooth, \(\mu\)-negligible hypersurface, which is the boundary relative to \(X\) of the open set \(\{f>t\}\), so that \(\Per_\mu(\{f>t\};X)=\int_{M_t}e^W\,\mathrm{d}\mathcal{H}^{N-1}\).
\end{lemma}

We now recall Talenti's theorem~\cite[Main Thm.]{Tal97}, see also~\cite[Thm.~3.10]{BS26}, which proves the P\'olya--Szeg\H{o} principle~\eqref{eq:PS} for bounded Lipschitz functions on \(X\).
When \(X\) is convex, it extends to Sobolev functions on an open set \(\Omega\subset X\) vanishing on the boundary of \(\Omega\) relative to \(X\), see~\cite[Thm.~5.1]{BS26}.

\begin{theorem}[P\'olya--Szeg\H{o} principle of~\cite{Tal97}]
\label{thm:talenti}
Assume~\ref{hp:N1}--\ref{hp:N3} and~\ref{hp:T1}--\ref{hp:T4}.
Then, for every bounded Lipschitz function \(u:X\to\R\) with \(\mu(\{u\neq0\})<+\infty\) and every \(p\in(1,+\infty)\),
\begin{equation}
\label{eq:PS}
\int_X|\nabla u^\sharp|^p\,\mathrm{d}\mu\le\int_X|\nabla u|^p\,\mathrm{d}\mu.
\end{equation}
\end{theorem}

Finally, we recall the structural result of~\cite[Thm.~4.1]{BS26}, which is the starting point of this paper.

\begin{theorem}
\label{thm:geo-uc}
Assume~\ref{hp:N1}--\ref{hp:N3}, \ref{hp:T3}, and~\ref{hp:T4}.
Then, the following hold.
\begin{enumerate}[label=(\alph*)]
\item\label{item:g}
There exists a continuous function \(g:(0,1)\to(0,+\infty)\) such that
\begin{equation}
\label{eq:g(f)}
\|\nabla f(x)\|=g(f(x)),
\qquad\forall\, x\in X_{\mathrm{reg}}.
\end{equation}
\item\label{item:separation}
Let \(c_0\in(0,1)\) and
\begin{equation}
\label{eq:G}
G(c)=\int_{c_0}^c\frac{\mathrm{d}t}{g(t)},
\qquad
c\in(0,1).
\end{equation}
If \(y_1,y_2\in X\) and \(0<\tau_1<\tau_2<1\) are such that \(f(y_1)\le\tau_1\) and \(f(y_2)\ge\tau_2\), then
\begin{equation}
\label{eq:separation}
\|y_1-y_2\| \ge G(\tau_2)-G(\tau_1) >0.
\end{equation}
\end{enumerate}
\end{theorem}

Whenever \cref{thm:geo-uc} applies, we fix \(c_0\in(0,1)\) once and for all, we let \(G\) be as in~\eqref{eq:G}, and we let
\begin{equation}
\label{eq:smin-smax}
s_{\min}=\lim_{c\to0^+}G(c)\in[-\infty,0),
\qquad
s_{\max}=\lim_{c\to1^-}G(c)\in(0,+\infty].
\end{equation}
As \(G\) is a strictly increasing \(C^1\) map of \((0,1)\) onto \((s_{\min},s_{\max})\), we extend its inverse to the whole \(\R\) by setting
\begin{equation}
\label{eq:G-inv}
G^{-1}(s)=0\quad\text{for }s\le s_{\min},
\qquad
G^{-1}(s)=1\quad\text{for }s\ge s_{\max}.
\end{equation}
Since the extended \(G^{-1}\) equals \(0\) on \((-\infty,s_{\min}]\), equals \(1\) on \([s_{\max},+\infty)\), and is strictly increasing on \((s_{\min},s_{\max})\), with \(G^{-1}(G(t))=t\), we have, for every \(t\in(0,1)\),
\begin{equation}
\label{eq:G-inv-sup}
G^{-1}(u)>t
\quad\Longleftrightarrow\quad
u>G(t),
\qquad
\forall\,u\in\R.
\end{equation}
In particular, if \(f=G^{-1}\circ h\) on \(X\) for some function \(h:X\to\R\), then \(\{f>t\}=\{h>G(t)\}\) for all \(t\in(0,1)\).

We conclude this section with an elementary fact on \(X_{\mathrm{reg}}\), which will be used throughout.

\begin{lemma}[Regular points]
\label{lem:reg}
Assume~\ref{hp:N1}.
Then, \(X_{\mathrm{reg}}\) is dense in \(X\).
If, in addition, \ref{hp:N2}, \ref{hp:N3}, \ref{hp:T3}, and~\ref{hp:T4} hold, then \(\{f=0\}=\emptyset\) if \(s_{\min}=-\infty\), and \(\{f=1\}=\emptyset\) if \(s_{\max}=+\infty\).
\end{lemma}

\begin{proof}
Let \(y\in X\) and \(r>0\) be such that \(B_r(y)\subset X\), and assume by contradiction that \(B_r(y)\cap X_{\mathrm{reg}}=\emptyset\).
Then, \(f\) only takes the values \(0\) and \(1\) on the connected open ball \(B_r(y)\), so that, by continuity, \(f\equiv f(y)\) on \(B_r(y)\).
Hence, \(\mu\bigl(f^{-1}(f(y))\bigr)\ge\mu(B_r(y))>0\), contradicting~\ref{hp:N1}.

As for the second part, let \(y_1\in X\) with \(f(y_1)=0\), and fix \(y_2\in X_{\mathrm{reg}}\), which exists by the first part.
For every \(\tau_1\in(0,f(y_2))\), \eqref{eq:separation} applied with \(\tau_2=f(y_2)\) gives \(\|y_1-y_2\|\ge G(f(y_2))-G(\tau_1)\).
Letting \(\tau_1\to0^+\), we get \(s_{\min}\ge G(f(y_2))-\|y_1-y_2\|>-\infty\).
The case \(f(y_1)=1\) is analogous.
\end{proof}

\section{Distance representation of \texorpdfstring{\(f\)}{f}}
\label{sec:structure}

If the lines of steepest ascent of \(f\) sweep the whole range of \(f\), that is, if~\ref{hp:SW} holds, then \(f\) can be determined by any of its level sets and by the profile \(g\) appearing in \cref{thm:geo-uc}\ref{item:g}, as the composition of \(G^{-1}\) with the signed distance function from a fixed level set. 
In particular, the level sets of \(f\) are parallel hypersurfaces.
By \cref{lem:sweep} below, this is always the case under the assumptions of Talenti's \cref{thm:talenti}.
In this section, \(c_0\), \(G\), \(s_{\min}\), \(s_{\max}\), and the extension of \(G^{-1}\) are as in~\eqref{eq:G}, \eqref{eq:smin-smax}, and~\eqref{eq:G-inv}.
We first show that the level sets of \(f\) are obtained by moving a given level set \(M_{c_0}\) along its normal; this will be used repeatedly in the sequel.

\begin{lemma}[Normal coordinates on \(X_{\mathrm{reg}}\)]
\label{lem:normal}
Assume~\ref{hp:N1}--\ref{hp:N3}, \ref{hp:T3}, \ref{hp:T4}, and~\ref{hp:SW}.
For \(x\in M_{c_0}\) and \(s\in(s_{\min},s_{\max})\), let \(P_s(x)=x+s\mathbf{N}(x)\).
Then, the following hold.
\begin{enumerate}[label=(\alph*)]
\item\label{item:flow-forward} For all \(x\in M_{c_0}\) and \(s\in(s_{\min},s_{\max})\), the point \(P_s(x)\) belongs to \(X_{\mathrm{reg}}\), \(\mathbf{N}(P_s(x))=\mathbf{N}(x)\), and
\begin{equation}
\label{eq:flow}
f(P_s(x))=G^{-1}(s).
\end{equation}
\item\label{item:flow-backward} Every \(y\in X_{\mathrm{reg}}\) can be uniquely written as \(y=P_s(x)\), with \(x\in M_{c_0}\) and \(s\in(s_{\min},s_{\max})\); moreover, \(s=G(f(y))\) and \(\mathbf{N}(y)=\mathbf{N}(x)\).
\end{enumerate}
In particular, the level sets of \(f\) are the parallel hypersurfaces \(M_{G^{-1}(s)}=P_s(M_{c_0})\), \(s\in(s_{\min},s_{\max})\).
\end{lemma}

\begin{proof}
Let \(\gamma\) be a maximal line of steepest ascent in \(X_{\mathrm{reg}}\), \emph{i.e.}, a maximal solution in \(X_{\mathrm{reg}}\) of \(\gamma'=\mathbf{N}(\gamma)\), which exists and is uniquely determined by \(\gamma(0)\), as \(\mathbf{N}\) is of class \(C^1\) on \(X_{\mathrm{reg}}\); let \(J\) be its interval of definition.
By~\ref{hp:T4}, \(\gamma''=0\), so that \(\mathbf{N}(\gamma)=\gamma'\) is constant on \(J\).
By~\eqref{eq:g(f)},
\[
\frac{\mathrm{d}}{\mathrm{d}s}f(\gamma(s))=\nabla f(\gamma(s))\cdot\mathbf{N}(\gamma(s))=\|\nabla f(\gamma(s))\|=g(f(\gamma(s))),
\]
so that \(\frac{\mathrm{d}}{\mathrm{d}s}G(f(\gamma(s)))=1\), that is,
\begin{equation}
\label{eq:G(f)}
G(f(\gamma(s)))=a+s,
\qquad\forall\, s\in J,
\qquad\text{where }a=G(f(\gamma(0))).
\end{equation}
By~\ref{hp:SW}, \(f(\gamma(J))=(0,1)\), and \(G\) maps \((0,1)\) onto \((s_{\min},s_{\max})\); hence, by~\eqref{eq:G(f)}, \(a+J=(s_{\min},s_{\max})\), that is,
\[
J=(s_{\min}-a,\,s_{\max}-a).
\]

\medskip
\textit{Step~1: Proof of~\ref{item:flow-forward}.}
If \(\gamma(0)=x\in M_{c_0}\), then \(a=G(c_0)=0\).
Therefore, \(J=(s_{\min},s_{\max})\) and \(\gamma(s)=x+s\mathbf{N}(x)=P_s(x)\) for all \(s\in J\).
Hence, \(P_s(x)\in X_{\mathrm{reg}}\), \(\mathbf{N}(P_s(x))=\mathbf{N}(x)\), as \(\mathbf{N}\) is constant along \(\gamma\), and~\eqref{eq:G(f)} gives~\eqref{eq:flow}.

\medskip
\textit{Step~2: Proof of~\ref{item:flow-backward}.}
Let \(y\in X_{\mathrm{reg}}\), and let \(\gamma\) be the maximal line of steepest ascent with \(\gamma(0)=y\), so that \(a=G(f(y))\in(s_{\min},s_{\max})\) and \(J\) contains \(-a\), as \(s_{\min}<0<s_{\max}\).
Let \(x=\gamma(-a)=y-a\mathbf{N}(y)\).
By~\eqref{eq:G(f)}, \(G(f(x))=0=G(c_0)\), so that \(f(x)=c_0\) by injectivity of \(G\), that is, \(x\in M_{c_0}\).
Moreover, as \(\mathbf{N}\) is constant along \(\gamma\), we have \(\mathbf{N}(x)=\mathbf{N}(y)\), whence \(y=x+a\mathbf{N}(x)=P_a(x)\).
As for uniqueness, if \(y=P_s(x)=P_{s'}(x')\), with \(x,x'\in M_{c_0}\) and \(s,s'\in(s_{\min},s_{\max})\), then \(s=s'=G(f(y))\) by~\eqref{eq:flow}, and \(\mathbf{N}(x)=\mathbf{N}(y)=\mathbf{N}(x')\) by~\ref{item:flow-forward}, whence \(x=y-s\mathbf{N}(y)=x'\).
\end{proof}

We can now prove the distance representation of \(f\).

\begin{theorem}
\label{thm:dist-rep}
Assume~\ref{hp:N1}--\ref{hp:N3}, \ref{hp:T3}, \ref{hp:T4}, and~\ref{hp:SW}, and let \(d\) be the signed distance from \(M_{c_0}\), positive in \(\{f>c_0\}\).
Then, the following hold.
\begin{enumerate}[label=(\alph*)]
\item\label{item:rep} The level sets of \(f\) are parallel hypersurfaces; specifically, one has
\[
f(y)=G^{-1}(d(y)),
\qquad
\forall\,y\in X.
\]
\item\label{item:focal} If \(s_{\min}>-\infty\), letting
\begin{equation}
\label{eq:F0}
F_0=\{\,x+s_{\min}\mathbf{N}(x)\,:\,x\in M_{c_0}\,\}\subset\overline{X},
\end{equation}
one has
\[
f(y)=G^{-1}(s_{\min}+\dist(y;F_0)),
\qquad
\forall\,y\in X.
\]
Analogously, if \(s_{\max}<+\infty\), letting
\begin{equation}
\label{eq:F1}
F_1=\{\,x+s_{\max}\mathbf{N}(x)\,:\,x\in M_{c_0}\,\}\subset\overline{X},
\end{equation}
one has \(f(y)=G^{-1}(s_{\max}-\dist(y;F_1))\) for all \(y\in X\).
\item\label{item:nonconstant} If \(\mathbf{N}\) is not constant on \(M_{c_0}\), then either \(s_{\min}>-\infty\) or \(s_{\max}<+\infty\).
\end{enumerate}
\end{theorem}

\begin{proof}
We divide the proof into three steps.

\medskip
\textit{Step~1: Proof of~\ref{item:rep}.}
Let \(y\in X_{\mathrm{reg}}\) and let \(s=G(f(y))\).
Since \(G\) is strictly increasing and \(G(c_0)=0\), the sign of \(s\) is that of \(f(y)-c_0\).
We claim that
\begin{equation}
\label{eq:dist-is-s}
\dist(y;M_{c_0})=|s|.
\end{equation}

On the one hand, by \cref{lem:normal}\ref{item:flow-backward}, there exists \(x_0\in M_{c_0}\) such that we have \(y=x_0+s\mathbf{N}(x_0)\); hence, \(\dist(y;M_{c_0})\le\|y-x_0\|=|s|\).

On the other hand, let \(x\in M_{c_0}\).
If \(s=0\), then trivially \(\|y-x\|\ge0=|s|\).
If \(s>0\), that is, \(c_0<f(y)\), we apply~\eqref{eq:separation} with
\[
y_1=x,\quad y_2=y,\quad \tau_1=c_0,\quad \tau_2=f(y),
\]
which is allowed, since \(0<\tau_1<\tau_2<1\), \(f(y_1)=\tau_1\), and \(f(y_2)=\tau_2\); this gives \(\|y-x\|\ge G(f(y))-G(c_0)=s\).
If \(s<0\), that is, \(f(y)<c_0\), we apply~\eqref{eq:separation} as above, but with the roles of \(x\) and \(y\) exchanged.
Thus
\[
\|y-x\|\ge |s|, \qquad \forall x\in M_{c_0}.
\]
Taking the infimum over \(x\in M_{c_0}\) proves the opposite inequality, giving~\eqref{eq:dist-is-s}.

Finally, by definition, \(d(y)=\dist(y;M_{c_0})\) if \(f(y)\ge c_0\), that is, if \(s\ge0\), and \(d(y)=-\dist(y;M_{c_0})\) otherwise.
Hence, by~\eqref{eq:dist-is-s}, \(d(y)=s=G(f(y))\), that is,
\begin{equation}
\label{eq:f-distance-rep}
f(y) = G^{-1}(d(y)),
\qquad \forall\, y \in X_{\mathrm{reg}}.
\end{equation}
The signed distance \(d\) is continuous on \(X\), and so is the extended function \(G^{-1}\) on \(\R\).
Hence, both sides of~\eqref{eq:f-distance-rep} are continuous on \(X\) and, coinciding on \(X_{\mathrm{reg}}\), which is dense in \(X\) by \cref{lem:reg}, they coincide on \(X\).

\medskip
\textit{Step~2: Proof of~\ref{item:focal}.}
Assume that \(s_{\min}>-\infty\), and let \(F_0\) be as in~\eqref{eq:F0}; notice that \(F_0\subset\overline{X}\), since \(P_s(x)\in X_{\mathrm{reg}}\) for all \(s\in(s_{\min},s_{\max})\) by \cref{lem:normal}\ref{item:flow-forward}.
Let \(y\in X_{\mathrm{reg}}\), and write
\[
y=x_y+s\mathbf{N}(x_y),
\quad
\text{with } x_y\in M_{c_0}
\text{ and }
s=G(f(y))=d(y),
\]
by \cref{lem:normal}\ref{item:flow-backward} and Step~1.

On the one hand, the point \(x_y+s_{\min}\mathbf{N}(x_y)\) belongs to \(F_0\) and has distance \(s-s_{\min}\) from \(y\), so that
\[
\dist(y; F_0) \le s-s_{\min}.
\]
On the other hand, let \(z=x+s_{\min}\mathbf{N}(x)\in F_0\), with \(x\in M_{c_0}\), and, for \(s_k\downarrow s_{\min}\), let \(z_k=P_{s_k}(x)\), so that \(z_k\to z\) and, by~\eqref{eq:flow}, \(f(z_k)=G^{-1}(s_k)\to0\).
Given \(\tau_1\in(0,f(y))\), for \(k\) large enough we have \(f(z_k)\le\tau_1\), and~\eqref{eq:separation}, applied with \(y_1=z_k\), \(y_2=y\), and \(\tau_2=f(y)\), yields
\[
\|y-z_k\|\ge G(f(y))-G(\tau_1)
=
s-G(\tau_1)
.
\]
Letting first \(k\to\infty\) and then \(\tau_1\to0^+\), we obtain \(\|y-z\|\ge s-s_{\min}\), and, passing to the infimum over \(z\in F_0\), we get \(\dist(y;F_0)\ge s-s_{\min}\).

This proves that \(\dist(y;F_0)=d(y)-s_{\min}\) and, by~\eqref{eq:f-distance-rep},
\[
f(y)=G^{-1}\bigl(s_{\min}+\dist(y;F_0)\bigr),
\qquad\forall\, y\in X_{\mathrm{reg}}.
\]
Since both sides are continuous and \(X_{\mathrm{reg}}\) is dense in \(X\) by \cref{lem:reg}, the equality holds on the whole \(X\).
The case \(s_{\max}<+\infty\) is analogous.

\medskip
\textit{Step~3: Proof of~\ref{item:nonconstant}.}
Assume that \(\mathbf{N}\) is not constant on \(M_{c_0}\).
Since \(M_{c_0}\) is connected, its shape operator \(S=-\mathrm{d}\mathbf{N}|_{TM_{c_0}}\) does not vanish identically; thus, there exist \(x_0\in M_{c_0}\) and an index \(j\) such that the \(j\)-th principal curvature \(\kappa_j(x_0)\) of \(M_{c_0}\) at \(x_0\), \emph{i.e.}, the \(j\)-th eigenvalue of \(S\), is not \(0\).
Along the line \(s\mapsto P_s(x_0)\), the principal curvatures (with respect to the arclength parameter) of the parallel hypersurfaces \(P_s(M_{c_0})=M_{G^{-1}(s)}\), see \cref{lem:normal}, satisfy the Riccati-type ODE
\[
\dot{\kappa}_i
=
\kappa_i^2,
\qquad
\kappa_i(0)=\kappa_i(x_0),
\]
see~\cite[Thm.~3]{Mag91}, whose solution is
\[
\kappa_i(x_0, s) = \frac{\kappa_i(x_0)}{1-s\kappa_i(x_0)}.
\]
Since the level sets \(M_{G^{-1}(s)}\) are smooth hypersurfaces for all \(s\in(s_{\min},s_{\max})\), their principal curvatures are finite, so that \(1-s\kappa_j(x_0)\neq0\) on \((s_{\min},s_{\max})\).
Hence, depending on the sign of \(\kappa_j(x_0)\),
\[
\text{either}\quad
s_{\min} \ge \frac{1}{\kappa_j(x_0)} >-\infty,
\qquad
\text{or}\quad
s_{\max} \le \frac{1}{\kappa_j(x_0)} < +\infty.
\qedhere
\]
\end{proof}

As a consequence, the superlevel sets of \(f\) are mutually parallel, and they fall into three families.
In the next statement, \(G\), \(s_{\min}\), and \(s_{\max}\) are as in~\eqref{eq:G} and~\eqref{eq:smin-smax}, and the focal sets \(F_0\) and \(F_1\) are as in~\eqref{eq:F0} and~\eqref{eq:F1}.

\begin{corollary}[Superlevel sets are parallel sets]
\label{cor:parallel}
Assume~\ref{hp:N1}--\ref{hp:N3}, \ref{hp:T3}, \ref{hp:T4}, and~\ref{hp:SW}, and let \(E_t=\{f>t\}\), for \(t\in(0,1)\).
Then,
\begin{equation}
\label{eq:parallel-sets}
E_{t'}=\bigl\{\,x\in X\,:\,\dist(x;X\setminus E_t)>G(t')-G(t)\,\bigr\},
\quad
0<t<t'<1,
\end{equation}
that is, each superlevel set of \(f\) is an inner parallel set, relative to \(X\), of any larger one.
Moreover, exactly one of the following holds.
\begin{enumerate}[label=(\alph*)]
\item\label{item:par-half} \(\mathbf{N}\) is constant on \(M_{c_0}\).
In this case, there exist \(\theta\in\mathbb{S}^{N-1}\) and \(r\in\R\) such that \(E_t=\{\,x\in X\,:\,\langle x,\theta\rangle>r+G(t)\,\}\) for all \(t\in(0,1)\).
\item\label{item:par-tubes} \(\mathbf{N}\) is not constant on \(M_{c_0}\), and \(s_{\max}<+\infty\).
In this case,
\[
E_t=\{\,x\in X\,:\,\dist(x;F_1)<s_{\max}-G(t)\,\},
\qquad\forall\,t\in(0,1).
\]
\item\label{item:par-complements} \(\mathbf{N}\) is not constant on \(M_{c_0}\), and \(s_{\max}=+\infty\).
In this case, \(s_{\min}>-\infty\), and
\[
E_t=\{\,x\in X\,:\,\dist(x;F_0)>G(t)-s_{\min}\,\},
\qquad\forall\,t\in(0,1).
\]
\end{enumerate}
\end{corollary}

\begin{proof}
\textit{Step~1: Parallel sets.}
Let \(0<t<t'<1\), and notice that the set \(X\setminus E_t=\{f\le t\}\) is nonempty by~\ref{hp:N1}.
We prove the two opposite inclusions.

Let \(y\in E_{t'}\), and pick \(\tau_2\in(t',1)\) with \(\tau_2\le f(y)\).
For every \(z\in X\setminus E_t\), that is, with \(f(z)\le t\), \eqref{eq:separation} applied with \(y_1=z\), \(y_2=y\), and \(\tau_1=t\) gives \(\|y-z\|\ge G(\tau_2)-G(t)\).
Hence, as \(G\) is strictly increasing,
\[
\dist(y;X\setminus E_t)\ge G(\tau_2)-G(t)>G(t')-G(t),
\]
that is, \(y\) belongs to the right-hand side of~\eqref{eq:parallel-sets}.

Conversely, let \(y\in X\setminus E_{t'}\), that is, \(f(y)\le t'\).
If \(f(y)\le t\), then \(y\in X\setminus E_t\), so that \(\dist(y;X\setminus E_t)=0\).
If \(t<f(y)\le t'\), then \(y\in X_{\mathrm{reg}}\) and, by \cref{lem:normal}\ref{item:flow-backward}, we can write \(y=P_s(x)\), where \(x\in M_{c_0}\) and \(s=G(f(y))\in(G(t),G(t')]\).
By~\eqref{eq:flow}, \(P_{G(t)}(x)\) belongs to \(M_t\subset X\setminus E_t\), and, as \(\|\mathbf{N}(x)\|=1\), its distance from \(y\) is \(s-G(t)\le G(t')-G(t)\).
In both cases, \(\dist(y;X\setminus E_t)\le G(t')-G(t)\), that is, \(y\) does not belong to the right-hand side of~\eqref{eq:parallel-sets}.

\medskip
\textit{Step~2: The three alternatives.}
The conditions on \(\mathbf{N}\) and \(s_{\max}\) in~\ref{item:par-half}--\ref{item:par-complements} are clearly mutually exclusive and exhaust all possibilities, so that it suffices to prove the corresponding descriptions of \(E_t\).
In each case, we show that \(f=G^{-1}\circ h\) on \(X\) for an explicit function \(h:X\to\R\), so that for all \(t\in(0,1)\) we have \(E_t=\{\,x\in X\,:\,h(x)>G(t)\,\}\) by~\eqref{eq:G-inv-sup}.

Proof of~\ref{item:par-half}: assume that \(\mathbf{N}\equiv\theta\) on \(M_{c_0}\).
By~\ref{hp:T3} and \cref{lem:Phi-f}, \(M_{c_0}\) is a connected smooth hypersurface, with constant normal \(\theta\); hence, \(M_{c_0}\subset\{\langle x,\theta\rangle=r\}\) for some \(r\in\R\).
By \cref{lem:normal}\ref{item:flow-backward} and \cref{thm:dist-rep}\ref{item:rep}, every \(y\in X_{\mathrm{reg}}\) can be written as \(y=x+d(y)\,\theta\), with \(x\in M_{c_0}\).
Therefore, we get \(d(y)=\langle y,\theta\rangle-r\).
Both sides being continuous on \(X\), and \(X_{\mathrm{reg}}\) being dense in \(X\) by \cref{lem:reg}, this identity holds on \(X\).
Hence, by \cref{thm:dist-rep}\ref{item:rep}, \(f=G^{-1}\circ h\) with \(h(x)=\langle x,\theta\rangle-r\), and~\ref{item:par-half} follows.

Proof of~\ref{item:par-tubes}: if \(\mathbf{N}\) is not constant on \(M_{c_0}\) and \(s_{\max}<+\infty\), then, by \cref{thm:dist-rep}\ref{item:focal}, \(f=G^{-1}\circ h\) with \(h=s_{\max}-\dist(\,\cdot\,;F_1)\), and~\ref{item:par-tubes} follows.

Proof of~\ref{item:par-complements}: if \(\mathbf{N}\) is not constant on \(M_{c_0}\) and \(s_{\max}=+\infty\), then \(s_{\min}>-\infty\) by \cref{thm:dist-rep}\ref{item:nonconstant}.
Then, owing to \cref{thm:dist-rep}\ref{item:focal}, we get \(f=G^{-1}\circ h\), where \(h=s_{\min}+\dist(\,\cdot\,;F_0)\), and~\ref{item:par-complements} follows.
\end{proof}

\begin{remark}[Talenti-symmetrized sets]
\label{rem:parallel}
Assume~\ref{hp:N1}--\ref{hp:N3}, and let \(E_t\) be as in \cref{cor:parallel}.
By \cref{lem:Phi-f}, the map \(t\mapsto\mu(E_t)=\Phi_f(t)\) is a continuous, strictly decreasing bijection of \((0,1)\) onto \((0,\mu(X))\), so that \(\{E_t\}_{t\in(0,1)}\) is a nested one-parameter family of open sets, in which every volume \(m\in(0,\mu(X))\) is attained by exactly one set \(E_t\).
By~\cite[Lem.~3.8]{BS26}, for every open set \(\Omega\subsetneq X\) with \(\mu(\Omega)<+\infty\), the Talenti-symmetrized set \(\Omega^\sharp\) is the member of this family with \(\mu(\Omega^\sharp)=\mu(\Omega)\), and, under~\ref{hp:T2}, all the sets \(E_t\) are isoperimetric.
By \cref{cor:parallel}, under~\ref{hp:T3}, \ref{hp:T4}, and~\ref{hp:SW}, Talenti symmetrization can only be performed with respect to nested families of mutually parallel sets.
For instance, concentric balls as in \cite[Sect.~6.1]{BS26} fall into~\ref{item:par-tubes}, with \(F_1=\{0\}\), while parallel half-spaces as in \cite[Sect.~6.2]{BS26} fall into~\ref{item:par-half}.
\end{remark}

\begin{remark}[Comparison with transnormal functions]
\label{rem:transnormal}
A smooth nonconstant function \(F\) on a given Riemannian manifold \(M\) is called \emph{transnormal} if \(\|\nabla F\|^2=b(F)\) for some smooth function \(b\), see~\cite{Bol73,Wan87}.
For such a function, the integral curves of \(\nabla F/\|\nabla F\|\) are geodesics and the regular level sets are parallel hypersurfaces; moreover, if \(M\) is complete, the focal varieties \(M_\pm\), \emph{i.e.}, the sets where \(F\) attains its minimum and its maximum, are smooth submanifolds, and each regular level set is a tube over either of them~\cite{Wan87}.

By \cref{thm:geo-uc}\ref{item:g}, the functions \(f\) considered here satisfy \(\|\nabla f\|^2=g(f)^2\) on \(X_{\mathrm{reg}}\), and \(g\) is smooth, as, along each line of steepest ascent \(\gamma\), the function \(s\mapsto f(\gamma(s))\) is a smooth diffeomorphism onto an open interval, and \(g(f(\gamma(s)))=\|\nabla f(\gamma(s))\|\).
Thus, \(f\) is transnormal on \(X_{\mathrm{reg}}\), with \(b=g^2\).
Conversely, transnormality implies that the lines of steepest ascent are straight.
Hence, under~\ref{hp:T3}, condition~\ref{hp:T4} is equivalent to the transnormality of \(f\) on \(X_{\mathrm{reg}}\), and \cref{thm:geo-uc}, \cref{lem:normal}, and \cref{thm:dist-rep}, as well as \cref{prop:focal-point}, can be regarded as counterparts of the above classical facts.
However, they do not follow from them, for the following reasons.

First, the function \(f\) is defined on a convex open set \(X\subset\R^N\) with boundary, rather than on a complete manifold.
A priori, its lines of steepest ascent may leave \(X\), which is why~\ref{hp:SW} is assumed in \cref{lem:normal,thm:dist-rep}.

Second, the focal sets \(F_0\) and \(F_1\) of \cref{thm:dist-rep} lie in \(\overline{X}\) and, unless \(X=\R^N\), they typically lie on \(\partial X\), where \(f\) is not defined.
Hence, the structure results of~\cite{Wan87} on the focal varieties do not apply, and the representation of \(f\) in terms of the distance from \(F_0\) or \(F_1\) has to be proved directly.
\end{remark}

\section{Isoperimetry and classification}
\label{sec:isoperimetry}

In this section, we combine the structure of \(f\) with the isoperimetric property of its superlevel sets.
We start by proving several preliminary properties that will lead to the proof of the \cref{thm:main}.

Throughout this section, whenever~\ref{hp:N1}--\ref{hp:N3}, \ref{hp:T3}, and~\ref{hp:T4} hold, we let \(c_0\), \(G\), \(s_{\min}\), and \(s_{\max}\) be as in~\eqref{eq:G} and~\eqref{eq:smin-smax}, and, for \(s\in(s_{\min},s_{\max})\),
\begin{equation}
\label{eq:phi-U-V}
\begin{split}
\varphi
&
=G\circ f\quad\text{on }X_{\mathrm{reg}},
\\
U_s
&
=\{f>G^{-1}(s)\},
\\
V(s)
&
=\mu(U_s)=\Phi_f(G^{-1}(s)).
\end{split}
\end{equation}
By~\ref{hp:N3} and \cref{lem:Phi-f}, \(V(s)\in(0,+\infty)\) for all \(s\in(s_{\min},s_{\max})\).

Under these assumptions, by \cref{thm:geo-uc}, \(\varphi\) is of class \(C^1\) on \(X_{\mathrm{reg}}\), whereas by \cref{thm:geo-uc}\ref{item:g} and~\eqref{eq:G}, its gradient satisfies
\begin{equation}
\label{eq:grad-phi}
\nabla\varphi
=
\frac{\nabla f}{g(f)}
=
\frac{\nabla f}{\|\nabla f\|}
=
\mathbf{N},
\end{equation}
and~\eqref{eq:separation} gives
\begin{equation}
\label{eq:phi-lip}
|\varphi(y_1)-\varphi(y_2)|\le\|y_1-y_2\|,
\qquad\forall\,y_1,y_2\in X_{\mathrm{reg}}.
\end{equation}
By \cref{lem:reg}, \(X_{\mathrm{reg}}\) is dense in \(X\).
Hence, by~\eqref{eq:phi-lip}, \(\varphi\) uniquely extends to a \(1\)-Lipschitz function on \(X\), still denoted by \(\varphi\), so that~\eqref{eq:phi-lip} holds for all \(y_1,y_2\in X\).
Further, since \(f\) is continuous, \(\varphi=s_{\min}\) on \(\{f=0\}\) and \(\varphi=s_{\max}\) on \(\{f=1\}\), these sets being empty if \(s_{\min}=-\infty\) and \(s_{\max}=+\infty\), respectively, by \cref{lem:reg}.
Moreover, if \(\gamma\) is a line of steepest ascent of \(f\), that is, an integral curve of \(\mathbf{N}\) in \(X_{\mathrm{reg}}\), which is parametrized by arclength as \(\|\mathbf{N}\|=1\), then, by~\eqref{eq:grad-phi},
\[
(\varphi\circ\gamma)'
=
\langle\nabla\varphi(\gamma),\mathbf{N}(\gamma)\rangle
=
1;
\]
that is, \(\varphi\) increases at unit speed along \(\gamma\).

\subsection{Comparison and sweep}
We first show that the function \(V\) in~\eqref{eq:phi-U-V}, defined on \((s_{\min},s_{\max})\) with \(s_{\min}\) and \(s_{\max}\) as in~\eqref{eq:smin-smax}, solves an ordinary differential equation driven by the isoperimetric profile \(I\) of~\ref{hp:T1}, and that it can be compared with any function satisfying the corresponding differential inequality.
This comparison principle will be used several times in the sequel.

\begin{lemma}[Comparison with the profile]
\label{lem:ode}
Assume that~\ref{hp:N1}--\ref{hp:N3} and~\ref{hp:T1}--\ref{hp:T4} hold.
Then, \(V\) is of class \(C^1\) on \((s_{\min},s_{\max})\), with \(V'=-I(V)\).
Moreover, \(V(s)\to0\) as \(s\to s_{\max}\), and \(V(s)\to\mu(X)\) as \(s\to s_{\min}\).
Further, the function
\begin{equation}
\label{eq:Lambda}
\Lambda(v)=\int_{V(0)}^{v}\frac{\mathrm{d}u}{I(u)},
\qquad
v\in(0,\mu(X)),
\end{equation}
is of class \(C^1\) and strictly increasing, and
\begin{equation}
\label{eq:Lambda-V}
\begin{gathered}
\Lambda(V(s))=-s,
\qquad
\forall\,s\in(s_{\min},s_{\max}),
\\
\Lambda(0^+)=-s_{\max},
\qquad
\Lambda(\mu(X)^-)=-s_{\min}.
\end{gathered}
\end{equation}
In particular, for all \(s\in(s_{\min},s_{\max})\),
\begin{equation}
\label{eq:width}
s_{\max}-s=\int_0^{V(s)}\frac{\mathrm{d}v}{I(v)},
\qquad
s-s_{\min}=\int_{V(s)}^{\mu(X)}\frac{\mathrm{d}v}{I(v)}.
\end{equation}
Finally, let \(J\subset\R\) be an interval, and let \(w\colon J\to(0,\mu(X))\) be absolutely continuous on every compact subinterval of \(J\).
Then, for all \(\tau_1,\tau_2\in J\) with \(\tau_1<\tau_2\),
\begin{equation}
\label{eq:comparison}
\begin{aligned}
\Lambda(w(\tau_2))&\le\Lambda(w(\tau_1))-(\tau_2-\tau_1),
&&\text{if }w'\le-I(w)\text{ a.e.\ in }J,
\\
\Lambda(w(\tau_2))&\ge\Lambda(w(\tau_1))+(\tau_2-\tau_1),
&&\text{if }w'\ge I(w)\text{ a.e.\ in }J.
\end{aligned}
\end{equation}
\end{lemma}

\begin{proof}
Since \(\varphi\) is Lipschitz, with \(\|\nabla\varphi\|=1\) on \(X_{\mathrm{reg}}\), for \(s_{\min}<s_1<s_2<s_{\max}\), by the coarea formula, see~\cite[Sect.~13.1]{Mag12book}, and by~\ref{hp:T2}, we have
\[
V(s_1)-V(s_2)
=
\mu(U_{s_1}\setminus U_{s_2})
=
\int_{s_1}^{s_2}\Per_\mu(U_\sigma;X)\,\mathrm{d}\sigma
=
\int_{s_1}^{s_2}I(V(\sigma))\,\mathrm{d}\sigma.
\]
Since \(I\) is continuous by~\ref{hp:T1}, the first claim follows.
The second one is a consequence of the definition of \(V\) and of \cref{lem:Phi-f}.

By~\ref{hp:T1}, \(\Lambda'=1/I\) is continuous and positive on \((0,\mu(X))\), so that \(\Lambda\) is of class \(C^1\) and strictly increasing.
As \(V\) takes values in \((0,\mu(X))\), we have \((\Lambda\circ V)'=V'/I(V)=-1\) on \((s_{\min},s_{\max})\).
Further, since \(0\in(s_{\min},s_{\max})\) and \(\Lambda(V(0))=0\), we have \(\Lambda(V(s))=-s\) for all \(s\in(s_{\min},s_{\max})\).
Letting \(s\to s_{\max}\) and \(s\to s_{\min}\), by the monotonicity of \(\Lambda\) and by the limits of \(V\), we get the last two identities in~\eqref{eq:Lambda-V}, and~\eqref{eq:width} follows by difference.

Finally, let \(w\) be as in the statement, and let \(\tau_1,\tau_2\in J\) with \(\tau_1<\tau_2\).
As \(w\) is continuous, on \([\tau_1,\tau_2]\) it takes values in a compact subset of \((0,\mu(X))\), on which \(\Lambda\) is Lipschitz.
Hence, \(\Lambda\circ w\) is absolutely continuous on \([\tau_1,\tau_2]\), with \((\Lambda\circ w)'=w'/I(w)\) a.e.
If \(w'\le-I(w)\) a.e., then \((\Lambda\circ w)'\le-1\) a.e., and the first inequality in~\eqref{eq:comparison} follows by integration on \([\tau_1,\tau_2]\); the second one is analogous.
\end{proof}

The next result shows that, as soon as the superlevel sets of \(f\) are isoperimetric, they are mutually parallel, and the lines of steepest ascent of \(f\) sweep its whole range.
In particular,~\ref{hp:SW} follows from the assumptions of \cref{thm:talenti}.
Therefore, under~\ref{hp:N1}--\ref{hp:N3} and~\ref{hp:T1}--\ref{hp:T4}, \cref{thm:dist-rep,cor:parallel} apply.
For the reader's convenience, we recall that \(\varphi\) and \(U_s\) are defined in~\eqref{eq:phi-U-V}, and \(s_{\min}\) and \(s_{\max}\) in~\eqref{eq:smin-smax}.

\begin{lemma}[Isoperimetry forces the sweep]
\label{lem:sweep}
Assume that~\ref{hp:N1}--\ref{hp:N3} and \ref{hp:T1}--\ref{hp:T4} hold.
Then, for all \(s\in(s_{\min},s_{\max})\) and \(x\in X_{\mathrm{reg}}\), we have
\begin{equation}
\label{eq:dist-identities}
\begin{split}
\dist(x;X\setminus U_s)&=\varphi(x)-s,\ \ \text{if }\varphi(x)>s,
\\
\dist(x;U_s)&=s-\varphi(x),\ \ \text{if }\varphi(x)<s.
\end{split}
\end{equation}
Moreover,~\ref{hp:SW} holds, that is, \(f\) takes all values in \((0,1)\) along each of its maximal lines of steepest ascent in \(X_{\mathrm{reg}}\).
\end{lemma}

\begin{proof}
Fix \(s\in(s_{\min},s_{\max})\).
Throughout the proof, we use that, as \(\varphi=s_{\min}\) on \(\{f=0\}\) and \(\varphi=s_{\max}\) on \(\{f=1\}\), the sets in~\eqref{eq:phi-U-V} can be written as
\begin{equation}
\label{eq:U-phi}
U_\sigma=\{\,x\in X\,:\,\varphi(x)>\sigma\,\},
\qquad
\forall\,\sigma\in(s_{\min},s_{\max}),
\end{equation}
which are relatively open in \(X\) as \(\varphi\) is continuous.

\medskip
\textit{Step~1: Inner parallel sets.}
Here we prove the first identity in~\eqref{eq:dist-identities}, by showing the two opposite inequalities.
Let \(\psi=\dist(\,\cdot\,;X\setminus U_s)\) on \(X\).

First, let \(x\in X\) and \(z\in X\setminus U_s\), so that \(\varphi(z)\le s\), by~\eqref{eq:U-phi}.
By the \(1\)-Lipschitz continuity of \(\varphi\), see~\eqref{eq:phi-lip}, we have \(\|x-z\|\ge\varphi(x)-\varphi(z)\ge\varphi(x)-s\), so that, taking the infimum over \(z\),
\begin{equation}
\label{eq:psi-lower}
\psi(x)\ge\varphi(x)-s,
\qquad
\forall\,x\in X.
\end{equation}

Second, we show that \(\psi(x)\le\varphi(x)-s\) for all \(x\in X_{\mathrm{reg}}\) with \(\varphi(x)>s\).
To this aim, let \(w(h)=\mu(\{\psi>h\})\) for \(h\ge0\).
As \(X\setminus U_s\) is relatively closed in \(X\), we have \(\{\psi>0\}=U_s\), so that it holds \(w(0)=V(s)\).
Moreover, by~\eqref{eq:psi-lower} and~\eqref{eq:U-phi}, we have \(U_{s+h}\subset\{\psi>h\}\) for all \(h\in(0,s_{\max}-s)\).
Since \(\psi\) is Lipschitz, with \(\|\nabla\psi\|=1\) a.e.\ in \(U_s\), by the coarea formula, for \(0\le h_1<h_2\), we have
\[
w(h_1)-w(h_2)
=
\int_{h_1}^{h_2}\Per_\mu(\{\psi>\sigma\};X)\,\mathrm{d}\sigma
\ge
\int_{h_1}^{h_2}I(w(\sigma))\,\mathrm{d}\sigma,
\]
where, in the last inequality, we used~\ref{hp:T1}, which applies as
\[
w(\sigma)
=
\mu(\{\psi>\sigma\})
\le
\mu(\{\psi>0\})
= V(s)<\mu(X).
\]
Hence, \(w\) is absolutely continuous, with
\[
w'(\sigma)
=
-\Per_\mu(\{\psi>\sigma\};X)
\le
-I(w(\sigma)),
\qquad
\text{for a.e.\ }\sigma>0.
\]
Let now \(h\in(0,s_{\max}-s)\) be such that \(w(h)>0\).
As \(w\) is nonincreasing, on \([0,h]\) it takes values in \([w(h),V(s)]\subset(0,\mu(X))\), so that, by~\eqref{eq:comparison} with \(\tau_1=0\) and \(\tau_2=h\), and by~\eqref{eq:Lambda-V}, with \(\Lambda\) as in~\eqref{eq:Lambda},
\[
\Lambda(w(h))
\le
\Lambda(V(s))-h
=
-(s+h)
=
\Lambda(V(s+h)).
\]
As \(\Lambda\) is strictly increasing, we get \(w(h)\le V(s+h)\) for all \(h\in(0,s_{\max}-s)\), which trivially holds also if \(w(h)=0\).
Since \(U_{s+h}\subset\{\psi>h\}\), this means that
\[
\mu\big(\{\psi>h\}\setminus U_{s+h}\big)
=
w(h)-V(s+h)
\le
0.
\]
As a consequence, the set \(\{\psi>h\}\setminus\overline{U_{s+h}}^X\), being open and \(\mu\)-negligible, is empty; in other words, by~\eqref{eq:U-phi},
\begin{equation}
\label{eq:inner-inclusion}
\psi(x)>h
\quad\Longrightarrow\quad
\varphi(x)\ge s+h,
\qquad
\forall\,x\in X,\ \forall\,h\in(0,s_{\max}-s).
\end{equation}
Let now \(x\in X_{\mathrm{reg}}\) with \(\varphi(x)>s\).
By~\eqref{eq:inner-inclusion}, \(\psi(x)\le h\) for all \(h\) in the interval \((\varphi(x)-s,s_{\max}-s)\), that is, \(\psi(x)\le\varphi(x)-s\), which, together with~\eqref{eq:psi-lower}, proves the first identity in~\eqref{eq:dist-identities}.

\medskip
\textit{Step~2: Outer parallel sets.}
The second identity in~\eqref{eq:dist-identities} follows along the lines of Step~1, with the sets \(\{\psi>h\}\subset U_s\) replaced by \(\{\psi<h\}\supset U_s\), where now \(\psi=\dist(\,\cdot\,;U_s)\) on \(X\); we only sketch the argument.
Analogously to Step~1, by the \(1\)-Lipschitz continuity of \(\varphi\), we have \(\psi\ge s-\varphi\) on \(X\), so that, by~\eqref{eq:U-phi}, \(\{\psi<h\}\subset U_{s-h}\) for all \(h\in(0,s-s_{\min})\).
Hence, letting now \(w(h)=\mu(\{\psi<h\})\), we have \(w(h)\le V(s-h)<\mu(X)\), so that~\ref{hp:T1} applies to \(\{\psi<h\}\), and the coarea formula gives that \(w\) is absolutely continuous on compact subintervals of \((0,s-s_{\min})\), with \(w'\ge I(w)\) a.e.
Finally, since \(\overline{U_s}^X\setminus U_s\subset M_{G^{-1}(s)}\) is \(\mu\)-negligible, we have \(w(0^+)=V(s)\), and the rest of the argument is as in Step~1, with reversed inequalities.

\medskip
\textit{Step~3: Sweep.}
Let \(x\in X_{\mathrm{reg}}\).
We first show that \(f\) takes all values in \([f(x),1)\) along the line of steepest ascent through \(x\).
Let \(s\in(\varphi(x),s_{\max})\).
By the second identity in~\eqref{eq:dist-identities}, there exists \(z\in\overline{U_s}\subset\overline{X}\) such that \(\|z-x\|=s-\varphi(x)\).
Let
\[
e=\frac{z-x}{\|z-x\|},
\qquad
y_\tau=x+\tau e,
\quad
\tau\in[0,s-\varphi(x)).
\]
Since \(X\) is convex and open, \(x\in X\), and \(z\in\overline{X}\), we have \(y_\tau\in X\).
We claim that
\begin{equation}
\label{eq:phi-along-segment}
\varphi(y_\tau)=\varphi(x)+\tau,
\qquad
\forall\,\tau\in[0,s-\varphi(x)).
\end{equation}
On the one hand, \(\varphi(y_\tau)\le\varphi(x)+\tau\) by~\eqref{eq:phi-lip}.
On the other hand, for every \(\xi\in U_s\), by~\eqref{eq:U-phi} and~\eqref{eq:phi-lip}, we have \(\|y_\tau-\xi\|\ge\varphi(\xi)-\varphi(y_\tau)>s-\varphi(y_\tau)\), so that, as \(z\in\overline{U_s}\),
\[
s-\varphi(y_\tau)
\le
\dist(y_\tau;U_s)
\le
\|y_\tau-z\|
=
s-\varphi(x)-\tau,
\]
where we used that \(y_\tau-z=(\tau-\|z-x\|)\,e\), that is, \(y_\tau\) lies on the segment from \(x\) to \(z\).
Rearranging, \(\varphi(y_\tau)\ge\varphi(x)+\tau\).
This proves~\eqref{eq:phi-along-segment}.
In particular, \(\varphi(y_\tau)\in[\varphi(x),s)\subset(s_{\min},s_{\max})\), so that \(y_\tau\in X_{\mathrm{reg}}\), as \(\varphi=s_{\min}\) on \(\{f=0\}\) and \(\varphi=s_{\max}\) on \(\{f=1\}\).
Since \(\|\nabla\varphi\|=1\) on \(X_{\mathrm{reg}}\), \eqref{eq:phi-along-segment} gives \(\nabla\varphi(y_\tau)=e\), that is, \(\mathbf{N}(y_\tau)=e\), for all \(\tau\in[0,s-\varphi(x))\).
Hence, the segment \(\tau\mapsto y_\tau\) is an integral curve of \(\mathbf{N}\), so that it is contained in the line of steepest ascent through \(x\), along which \(\varphi\) therefore takes all values in \([\varphi(x),s)\).
As \(s<s_{\max}\) is arbitrary, \(\varphi\) takes all values in \([\varphi(x),s_{\max})\), that is, \(f\) takes all values in \([f(x),1)\), along this line.

Analogously, letting \(s\in(s_{\min},\varphi(x))\) and exploiting the first identity in~\eqref{eq:dist-identities}, one shows that \(f\) takes all values in \((0,f(x)]\) along the same line.
This proves~\ref{hp:SW}.
\end{proof}

We can now prove, one at a time, the remaining ingredients of the \cref{thm:main}.

\subsection{Balls}
First, we show that, if \(s_{\max}<+\infty\), defined in~\eqref{eq:smin-smax}, the superlevel sets of \(f\) are concentric balls; in particular, the focal set \(F_1\) in~\eqref{eq:F1} reduces to a point.
Here, \(G\) is as in~\eqref{eq:G}.

\begin{proposition}[The focal set \(F_1\) is a point]
\label{prop:focal-point}
Assume that~\ref{hp:N1}--\ref{hp:N3} and \mbox{\ref{hp:T1}--\ref{hp:T4}} hold, and that \(s_{\max}<+\infty\).
Then, there exists \(x_0\in\overline{X}\) such that
\[
\{\,x\in X\,:\,f(x)>t\,\}=B_{s_{\max}-G(t)}(x_0)\cap X,
\qquad\forall\,t\in(0,1),
\]
that is, the set \(F_1\) in~\eqref{eq:F1} reduces to \(\{x_0\}\); in particular, \(\mathbf{N}\) is not constant on \(M_{c_0}\).
\end{proposition}

\begin{proof}
By \cref{lem:sweep},~\ref{hp:SW} holds, so that \cref{thm:dist-rep} applies.
By \cref{thm:dist-rep}\ref{item:focal}, we have \(F_1\subset\overline{X}\) and \(f=G^{-1}(s_{\max}-\dist(\,\cdot\,;F_1))\) on \(X\); hence, by~\eqref{eq:phi-U-V} and~\eqref{eq:G-inv-sup}, the superlevel sets of \(f\) are tubes around \(F_1\), that is,
\begin{equation}
\label{eq:tubes}
U_{s_{\max}-r}=\{\,x\in X\,:\,\dist(x;F_1)<r\,\},
\qquad
\forall\,r\in(0,s_{\max}-s_{\min}).
\end{equation}
%

\medskip
\textit{Step~1: Balls grow faster than tubes.}
Here we prove that
\begin{equation}
\label{eq:ball-vs-tube}
\mu(B_r(a)\cap X)\ge V(s_{\max}-r),
\quad\forall\,a\in\overline{X},\ \forall\, r\in(0,s_{\max}-s_{\min}).
\end{equation}

Fix \(a\in\overline{X}\) and \(r\in(0,s_{\max}-s_{\min})\), let \(w(\sigma)=\mu(B_\sigma(a)\cap X)\) for \(\sigma\in[0,r]\), and assume by contradiction that \(w(r)<V(s_{\max}-r)\).
Then, \(0<w(\sigma)<\mu(X)\) for all \(\sigma\in(0,r]\): the lower bound holds as \(a\in\overline{X}\), and the upper one as \(w(\sigma)\le w(r)<V(s_{\max}-r)<\mu(X)\), owing to the definition of \(V\) in~\eqref{eq:phi-U-V}.
By the coarea formula,
\[
w(\sigma)
=
\int_0^\sigma\Per_\mu(B_\tau(a)\cap X;X)\,\mathrm{d}\tau,
\]
so that \(w\) is absolutely continuous on \([0,r]\) and, by~\ref{hp:T1}, \(w'\ge I(w)\) a.e.\ on \((0,r)\).
Hence, by~\eqref{eq:comparison} on \(J=(0,r]\), for every \(\sigma_0\in(0,r)\), recalling also that \(\Lambda\) is strictly increasing by \cref{lem:ode}, and by~\eqref{eq:Lambda-V}, we have
\[
\Lambda(w(r))
\ge
\Lambda(w(\sigma_0))+r-\sigma_0
\ge
\Lambda(0^+)+r-\sigma_0
=
r-s_{\max}-\sigma_0.
\]
Letting \(\sigma_0\to0^+\), we get \(\Lambda(w(r))\ge r-s_{\max}=\Lambda(V(s_{\max}-r))\), the last equality owing to~\eqref{eq:Lambda-V}.
Since \(\Lambda\) is strictly increasing, it follows that \(w(r)\ge V(s_{\max}-r)\), a contradiction.
This proves~\eqref{eq:ball-vs-tube}.

\medskip
\textit{Step~2: Conclusion.}
Assume by contradiction that there exist \(a,b\in F_1\) with \(a\neq b\), and let \(r\in(0,s_{\max}-s_{\min})\) be such that \(2r<\|a-b\|\).
The sets \(B_r(a)\cap X\) and \(B_r(b)\cap X\) are disjoint and, by~\eqref{eq:tubes}, they are contained in \(U_{s_{\max}-r}\).
Hence, by~\eqref{eq:ball-vs-tube}, which applies as \(a,b\in F_1\subset\overline{X}\),
\[
V(s_{\max}-r)
\ge
\mu(B_r(a)\cap X)+\mu(B_r(b)\cap X)
\ge
2V(s_{\max}-r),
\]
which contradicts \(0<V(s_{\max}-r)<+\infty\).
Therefore, \(F_1=\{x_0\}\) for some \(x_0\in\overline{X}\), and, by~\eqref{eq:tubes} with \(r=s_{\max}-G(t)\), we get
\[
\{f>t\}=B_{s_{\max}-G(t)}(x_0)\cap X,
\qquad
\forall\,t\in(0,1).
\]
Finally, if \(\mathbf{N}\equiv\theta\) were constant on \(M_{c_0}\), then \(F_1=\{\,x+s_{\max}\theta\,:\,x\in M_{c_0}\,\}\) would be a translate of the hypersurface \(M_{c_0}\), and not a single point.
Hence, \(\mathbf{N}\) is not constant on \(M_{c_0}\).
\end{proof}

\subsection{Finite measure forces half-spaces}
We now show that the total measure \(\mu(X)\) restricts the alternatives of \cref{cor:parallel}: if \(\mu(X)<+\infty\), there are no focal sets, and only parallel half-spaces may occur.

\begin{proposition}[Finite total measure]
\label{prop:width}
Assume that~\ref{hp:N1}--\ref{hp:N3} and~\ref{hp:T1}--\ref{hp:T4} hold, and that \(\mu(X)<+\infty\).
Then, \(s_{\min}=-\infty\) and \(s_{\max}=+\infty\); in particular, \(\mathbf{N}\) is constant on \(M_{c_0}\), that is, the superlevel sets of \(f\) are parallel half-spaces intersected with \(X\).
\end{proposition}

\begin{proof}
By \cref{lem:sweep},~\ref{hp:SW} holds, so that \cref{thm:dist-rep} applies.

\medskip
\textit{Step~1: Symmetry of the profile.}
First, \(\Per_\mu(E;X)=\Per_\mu(X\setminus E;X)\) for every Borel set \(E\subset X\).
Second, for all \(t\in(0,1)\), the set \(X\setminus\{f>t\}\) has measure \(\mu(X)-\Phi_f(t)<\mu(X)\), as \(\mu(X)<+\infty\) and \(\Phi_f(t)>0\) by \cref{lem:Phi-f}.
Then, by~\ref{hp:T1} applied to \(X\setminus\{f>t\}\) and by~\ref{hp:T2}, we have
\[
I(\mu(X)-\Phi_f(t))\le\Per_\mu(\{f>t\};X)=I(\Phi_f(t)),
\qquad\forall\,t\in(0,1).
\]
As \(\Phi_f\) maps \((0,1)\) onto \((0,\mu(X))\) by \cref{lem:Phi-f}, we get \(I(\mu(X)-m)\le I(m)\) for all \(m\in(0,\mu(X))\), and, exchanging the roles of \(m\) and \(\mu(X)-m\),
\begin{equation}
\label{eq:I-symmetric}
I(m)=I(\mu(X)-m),
\qquad
\forall\,m\in(0,\mu(X)).
\end{equation}
By~\eqref{eq:width} and~\eqref{eq:I-symmetric}, \(s_{\max}<+\infty\) if and only if \(s_{\min}>-\infty\), both conditions being equivalent to the integrability of \(1/I\) near \(0\). 

\medskip
\textit{Step~2: No focal sets.}
We claim that \(s_{\max}=+\infty\), so that also \(s_{\min}=-\infty\), by Step~1.
This suffices: then, \cref{thm:dist-rep}\ref{item:nonconstant} implies that \(\mathbf{N}\) is constant on \(M_{c_0}\), and the conclusion follows from \cref{cor:parallel}\ref{item:par-half}.

Assume by contradiction that \(s_{\max}<+\infty\), so that \(s_{\min}>-\infty\) as well.
We show that \(\tilde f=1-f\) satisfies the assumptions of \cref{prop:focal-point}.
Indeed, \(\tilde f\) is smooth and satisfies~\ref{hp:N1}--\ref{hp:N3}, the latter since \(\mu(X)<+\infty\).
It has the same level sets as \(f\), which are connected, and the same lines of steepest ascent, with reversed orientation, so that~\ref{hp:T3} and~\ref{hp:T4} hold.
Moreover, for all \(t\in(0,1)\), as the level sets of \(f\) are \(\mu\)-negligible, we have \(\Phi_{\tilde f}(t)=\mu(X)-\Phi_f(1-t)\), and, as the \(\mu\)-perimeter is invariant under complementation in \(X\) and under modifications on \(\mu\)-negligible sets, by~\ref{hp:T2} and~\eqref{eq:I-symmetric}, we get
\begin{align*}
\Per_\mu(\{\tilde f>t\};X)
&
=
\Per_\mu(\{f<1-t\};X)
=
\Per_\mu(\{f\ge1-t\};X)
\\
&
=
\Per_\mu(\{f>1-t\};X)
=
I(\Phi_f(1-t))
\\
&
=
I(\mu(X)-\Phi_{\tilde f}(t))
=
I(\Phi_{\tilde f}(t)).
\end{align*}
Hence,~\ref{hp:T2} holds for \(\tilde f\).
Finally, \(\|\nabla\tilde f\|=\tilde g(\tilde f)\), with \(\tilde g(t)=g(1-t)\), so that, choosing \(1-c_0\) as reference level in~\eqref{eq:G}, the function associated with \(\tilde f\) is \(\tilde G(c)=-G(1-c)\), and the quantity \(s_{\max}\) associated with \(\tilde f\) is \(-s_{\min}<+\infty\).
Therefore, \cref{prop:focal-point} applies to \(\tilde f\).

By \cref{prop:focal-point}, applied to \(f\) and to \(\tilde f\), there exist \(x_0,x_1\in\overline{X}\) such that, for all \(t\in(0,1)\),
\begin{align*}
\{f>t\}
&
=B_{s_{\max}-G(t)}(x_1)\cap X,
\\
\{f<t\}
=
\{\tilde f>1-t\}
&
=
B_{G(t)-s_{\min}}(x_0)\cap X.
\end{align*}
By \cref{lem:Phi-f}, applied to \(f\) and to \(\tilde f\), the nonempty level set \(M_t\) is the boundary relative to \(X\) of both sets, that is,
\[
M_t
=
\partial B_{s_{\max}-G(t)}(x_1)\cap X
=
\partial B_{G(t)-s_{\min}}(x_0)\cap X,
\]
which is therefore a nonempty relatively open subset of both spheres.
Since a sphere containing a nonempty relatively open subset of another sphere coincides with it, we get \(x_0=x_1\) and \(s_{\max}-G(t)=G(t)-s_{\min}\) for all \(t\in(0,1)\), which is impossible, \(G\) being strictly increasing.
This contradiction shows that \(s_{\max}=+\infty\), which, as observed above, concludes the proof.
\end{proof}

\subsection{Complements of metric neighborhoods}
Next, we show that, whenever \(s_{\min}>-\infty\) and \(s_{\max}=+\infty\), as in alternative~\ref{item:par-complements} of \cref{cor:parallel}, the density is not integrable near the focal set \(F_0\).
Here, \(s_{\min}\) and \(s_{\max}\) are as in~\eqref{eq:smin-smax}, and the focal set \(F_0\) is as in~\eqref{eq:F0}.

\begin{proposition}[Non-integrability near \(F_0\)]
\label{prop:non-integrable}
Assume that~\ref{hp:N1}--\ref{hp:N3} and \mbox{\ref{hp:T1}--\ref{hp:T4}} hold, and that \(s_{\min}>-\infty\) and \(s_{\max}=+\infty\).
Then, \(\mu(X)=+\infty\) and
\[
\mu(\{\,x\in X\,:\,\dist(x;F_0)<r\,\})=+\infty,
\qquad\forall\, r>0.
\]
Further, \(\mu(B_\rho(a)\cap X)=+\infty\) for all \(a\in\overline{F_0}\) and all \(\rho>0\).
In particular, if there exists \(a\in F_0\) such that \(\mu(X\setminus B_\sigma(a))<+\infty\) for all \(\sigma>0\), then \(F_0=\{a\}\).
\end{proposition}

\begin{proof}
By \cref{lem:sweep},~\ref{hp:SW} holds, so that \cref{thm:dist-rep} applies; in particular, \(F_0\neq\emptyset\) and \(f=G^{-1}(s_{\min}+\dist(\,\cdot\,;F_0))\) on \(X\), by \cref{thm:dist-rep}\ref{item:focal}.

Further, since \(s_{\min}>-\infty\), by \cref{prop:width}, we have \(\mu(X)=+\infty\).
For every \(r>0\), as \(s_{\min}<s_{\min}+r<s_{\max}=+\infty\), the value \(t=G^{-1}(s_{\min}+r)\) lies in \((0,1)\) and satisfies \(G(t)=s_{\min}+r\), so that, by~\eqref{eq:G-inv-sup} and the discussion thereafter,
\begin{equation}
\label{eq:meas-dist-fG-1}
\{\,x\in X\,:\,\dist(x;F_0)>r\,\}=\{f>G^{-1}(s_{\min}+r)\},
\qquad\forall\,r>0.
\end{equation}

\medskip
\textit{Step~1: Tubes around \(F_0\) have infinite measure.}
Let \(r>0\).
As \(\mu(X)=+\infty\) and the right-hand side of~\eqref{eq:meas-dist-fG-1} has finite measure by~\ref{hp:N3}, then its complement in \(X\), which is contained in \(\{\dist(\,\cdot\,;F_0)<2r\}\), has infinite measure.

\medskip
\textit{Step~2: Balls centered on \(\overline{F_0}\) have infinite measure.}
Let \(a\in\overline{F_0}\), and assume by contradiction that \(\mu(B_\rho(a)\cap X)<+\infty\) for some \(\rho>0\).
If \(\overline{F_0}\subset B_{\rho/2}(a)\), then 
\[
\{\dist(\,\cdot\,;F_0)<\rho/2\}\cap X\subset B_\rho(a)\cap X
\]
has finite measure, against Step~1.
Therefore, \(C=\overline{F_0}\setminus B_{\rho/2}(a)\) is nonempty and closed, and we let
\[
w(\sigma)=\mu(\{\,x\in X\,:\,\dist(x;C)>\sigma\,\}),
\qquad
\sigma>0.
\]
By~\eqref{eq:meas-dist-fG-1}, the definition of \(V\) in~\eqref{eq:phi-U-V}, and \(s_{\max}=+\infty\), we have
\begin{align*}
\mu(\{\dist(\,\cdot\,;F_0)>\sigma\})
&
=
\mu(\{f>G^{-1}(s_{\min}+\sigma)\})
\\
&
=
\mu(U_{s_{\min}+\sigma})
=
V(s_{\min}+\sigma)
\in (0,+\infty),
\qquad
\forall \sigma>0.
\end{align*}

On the one hand, as \(C\subset\overline{F_0}\), we have \(\dist(\,\cdot\,;C)\ge\dist(\,\cdot\,;F_0)\), and thus \(w(\sigma)\ge V(s_{\min}+\sigma)>0\).
On the other hand, we claim that for all \(\sigma\in(0,\rho/4)\)
\begin{equation}
\label{eq:KinF0}
\{\,x\in X:\dist(x;C)>\sigma\,\}
\subset
\{\,x\in X:\dist(x;F_0)>\sigma\,\}
\cup
\bigl(B_\rho(a)\cap X\bigr).
\end{equation}
Indeed, let \(x\in X\) be such that \(\dist(x;C)>\sigma\); if \(\dist(x;F_0)>\sigma\), there is nothing to prove, so we may assume that \(\dist(x;F_0)\le\sigma\).
As \(\overline{F_0}\) is closed, there exists \(b\in\overline{F_0}\) such that \(\|x-b\|=\dist(x;F_0)\le\sigma\).
Further, \(b\notin C\), as \(\dist(x;C)>\sigma\); hence, \(b\in B_{\rho/2}(a)\), by definition of \(C\).
Therefore,
\[
\|x-a\|\le\|x-b\|+\|b-a\|<\sigma+\rho/2<\rho,
\]
that is, \(x\in B_\rho(a)\), which proves~\eqref{eq:KinF0}.
By~\eqref{eq:KinF0} and our contradiction assumption, we get
\[
w(\sigma)\le V(s_{\min}+\sigma)+\mu(B_\rho(a)\cap X)<+\infty.
\]

We now show that, in fact, \(w(\sigma)\le V(s_{\min}+\sigma)\) for all \(\sigma\in(0,\rho/4)\).
By the coarea formula applied to the Lipschitz function \(\dist(\,\cdot\,;C)\), whose gradient has unit length a.e.\ outside \(C\), and by~\ref{hp:T1}, which applies since \(w<+\infty=\mu(X)\) on \((0,\rho/4)\), for \(0<\sigma_0<\sigma<\rho/4\), we have
\[
w(\sigma_0)-w(\sigma)
=
\int_{\sigma_0}^{\sigma}\Per_\mu(\{\dist(\,\cdot\,;C)>\tau\};X)\,\mathrm{d}\tau
\ge
\int_{\sigma_0}^{\sigma}I(w(\tau))\,\mathrm{d}\tau,
\]
so that \(w\) is absolutely continuous on \([\sigma_0,\sigma]\), with \(w'\le-I(w)\) a.e.
As \(w\) takes values in \((0,+\infty)=(0,\mu(X))\), by~\eqref{eq:comparison} on \(J=(0,\rho/4)\), since the function \(\Lambda\) of \cref{lem:ode} is strictly increasing, and by~\eqref{eq:Lambda-V}, we get
\begin{align*}
\Lambda(w(\sigma))
&
\le
\Lambda(w(\sigma_0))-(\sigma-\sigma_0)
\\
&
\le
\Lambda(\mu(X)^-)-\sigma+\sigma_0
=
-s_{\min}-\sigma+\sigma_0.
\end{align*}
Letting now \(\sigma_0\to0^+\) and owing once again to~\eqref{eq:Lambda-V}, we finally obtain that \(\Lambda(w(\sigma))\le-(s_{\min}+\sigma)=\Lambda(V(s_{\min}+\sigma))\), that is, \(w(\sigma)\le V(s_{\min}+\sigma)\), as \(\Lambda\) is strictly increasing.

Hence, for \(\sigma\in(0,\rho/4)\), we have \(w(\sigma)=V(s_{\min}+\sigma)<+\infty\), so that the set \(\{\,x\in X\,:\,\dist(x;F_0)\le\sigma<\dist(x;C)\,\}\), whose measure is \(w(\sigma)-V(s_{\min}+\sigma)\), is \(\mu\)-negligible, and its open subset \(\{\,x\in X\,:\,\dist(x;F_0)<\sigma<\dist(x;C)\,\}\) is empty.
However, the latter contains \(B_\sigma(a)\cap X\), which is nonempty as \(a\in\overline{F_0}\subset\overline{X}\).
Indeed, let \(x\in B_\sigma(a)\cap X\).
On the one hand, as \(a\in\overline{F_0}\), we have \(\dist(x;F_0)\le\|x-a\|<\sigma\).
On the other hand, every \(c\in C\) satisfies \(\|c-a\|\ge\rho/2\), so that \(\|x-c\|\ge\|c-a\|-\|x-a\|>\rho/2-\sigma\); hence, \(\dist(x;C)\ge\rho/2-\sigma>\sigma\), as \(\sigma<\rho/4\).
This contradiction proves the second claim.

\medskip
\textit{Step~3: A criterion for \(F_0\) to be a point.}
Finally, let \(a\in F_0\) be such that \(\mu(X\setminus B_\sigma(a))<+\infty\) for all \(\sigma>0\), and assume by contradiction that there exists \(b\in\overline{F_0}\) with \(b\neq a\).
Then, for \(\rho=\|a-b\|/2\), the set \(B_\rho(b)\cap X\subset X\setminus B_\rho(a)\) has finite measure, against the second claim.
Hence, \(F_0=\{a\}\).
\end{proof}

If \(s_{\min}>-\infty\) and \(s_{\max}=+\infty\), the lines of steepest ascent of \(f\) are half-lines emanating from \(F_0\) and contained in \(X\).
Together with the convexity of \(X\), this forces the closure of the focal set to be convex and, combined with the smoothness of the level sets of \(f\), to be flat; moreover, its affine hull has codimension at least two unless \(\mathbf{N}\) is constant on \(M_{c_0}\).

\begin{proposition}[Focal sets are flat]
\label{prop:focal-c}
Assume that~\ref{hp:N1}--\ref{hp:N3} and~\ref{hp:T1}--\ref{hp:T4} hold, and that \(s_{\min}>-\infty\) and \(s_{\max}=+\infty\).
Let \(F_0\) be the focal set in~\eqref{eq:F0}, and let \(L\) be the affine hull of \(\overline{F_0}\).
Then, \(L\cap X=\emptyset\), \(\overline{F_0}=L\cap\overline{X}\), and \(\dist(\,\cdot\,;F_0)=\dist(\,\cdot\,;L)\) on \(X\); in particular,
\[
\{\,x\in X\,:\,f(x)>t\,\}=\{\,x\in X\,:\,\dist(x;L)>G(t)-s_{\min}\,\},
\quad\forall\,t\in(0,1).
\]
Moreover, \(L\) is a hyperplane if \(\mathbf{N}\) is constant on \(M_{c_0}\), and \(\dim L\le N-2\) otherwise.
\end{proposition}

\begin{proof}
By \cref{lem:sweep},~\ref{hp:SW} holds, so that \cref{lem:normal,thm:dist-rep} apply.
As \(s_{\min}>-\infty\), by \cref{thm:dist-rep}\ref{item:focal}, we have \(F_0\subset\overline{X}\), hence \(\overline{F_0}\subset\overline{X}\), and \(f=G^{-1}(s_{\min}+d_0)\) on \(X\), where \(d_0=\dist(\,\cdot\,;F_0)=\dist(\,\cdot\,;\overline{F_0})\).
In particular, by the discussion following~\eqref{eq:G-inv-sup}, \(\{f>t\}=\{d_0>G(t)-s_{\min}\}\) for all \(t\in(0,1)\), so that the formula for the superlevel sets of \(f\) follows as soon as we show that \(d_0=\dist(\,\cdot\,;L)\) on \(X\).

\medskip
\textit{Step~1: Rays from the focal set.}
By the above representation of \(f\), recalling that \(\varphi=G\circ f\), we have \(\varphi=s_{\min}+d_0\) on \(X_{\mathrm{reg}}\).
Hence, as \(s_{\max}=+\infty\), by \cref{lem:normal}\ref{item:flow-forward} and by~\eqref{eq:flow},
\begin{equation}
\label{eq:distance-from-f0}
P_s(x) \in X_{\mathrm{reg}}
\quad
\text{and}
\quad
d_0(P_s(x))=s-s_{\min},
\qquad
\forall s>s_{\min},
\, 
x\in M_{c_0}.
\end{equation}
Let now \(y\in X_{\mathrm{reg}}\).
By \cref{lem:normal}\ref{item:flow-backward}, we have \(y=P_s(x)=x+s\mathbf{N}(x)\) for some \(x\in M_{c_0}\), with \(s=\varphi(y)\) and \(\mathbf{N}(x)=\mathbf{N}(y)\).
We let \(\pi_0(y)=x+s_{\min}\mathbf{N}(x)\), which belongs to \(F_0\) by~\eqref{eq:F0}, so that, as \(\mathbf{N}(x)=\mathbf{N}(y)\),
\begin{equation}
\label{eq:flow-from-F0}
\pi_0(y)+\tau\mathbf{N}(y)=x+(s_{\min}+\tau)\mathbf{N}(x)=P_{s_{\min}+\tau}(x),
\qquad\forall\,\tau>0.
\end{equation}
Therefore, by~\eqref{eq:distance-from-f0}, the whole open half-line \(\{\,\pi_0(y)+\tau\mathbf{N}(y)\,:\,\tau>0\,\}\) is contained in \(X_{\mathrm{reg}}\), and
\begin{equation}
\label{eq:d0-on-ray}
d_0(\pi_0(y)+\tau\mathbf{N}(y))=\tau,
\qquad\forall\,\tau>0.
\end{equation}
For \(\tau=s-s_{\min}=d_0(y)\), the point in~\eqref{eq:flow-from-F0} is \(P_s(x)=y\), so that
\begin{equation}
\label{eq:y-pi0}
y=\pi_0(y)+d_0(y)\,\mathbf{N}(y).
\end{equation}

\medskip
\textit{Step~2: \(\overline{F_0}\) is convex.}
Fix \(y\in X_{\mathrm{reg}}\).
By~\eqref{eq:d0-on-ray}, for every \(\tau>0\), the open ball \(B_\tau(\pi_0(y)+\tau\mathbf{N}(y))\) does not meet \(\overline{F_0}\), and, since the union of these balls over \(\tau>0\) is the open half-space \(\{\,z\,:\,\langle z-\pi_0(y),\mathbf{N}(y)\rangle>0\,\}\), we get
\begin{equation}
\label{eq:support-half-spaces}
\overline{F_0}\subset\{\,z\in\R^N\,:\,\langle z-\pi_0(y),\mathbf{N}(y)\rangle\le0\,\}.
\end{equation}
This half-space being convex, it also contains the convex hull \(\co(\overline{F_0})\) of \(\overline{F_0}\), so that, for all \(z\in\co(\overline{F_0})\), by~\eqref{eq:y-pi0} and~\eqref{eq:support-half-spaces}, we have
\[
\|y-z\|
\ge
\langle y-z,\mathbf{N}(y)\rangle
=
d_0(y)+\langle\pi_0(y)-z,\mathbf{N}(y)\rangle
\ge
d_0(y).
\]
Passing to the infimum over \(z\in \co(\overline{F_0})\), we have \(\dist(y;\co(\overline{F_0}))\ge d_0(y)\); owing to the fact that the opposite inequality trivially holds, \(\dist(y;\co(\overline{F_0}))=d_0(y)\) for all \(y\in X_{\mathrm{reg}}\) and, by continuity and by the density of \cref{lem:reg}, for all \(y\in\overline{X}\).
As \(\overline{X}\) is convex and contains \(\overline{F_0}\), it contains \(\co(\overline{F_0})\), so that \(d_0(z)=\dist(z;\co(\overline{F_0}))=0\) for all \(z\in\co(\overline{F_0})\).
Hence, \(\co(\overline{F_0})\subset\overline{F_0}\), that is, \(\overline{F_0}\) is convex.

\medskip
\textit{Step~3: \(\overline{F_0}\subset\partial X\).}
By \cref{prop:non-integrable}, \(\mu(B_\rho(a)\cap X)=+\infty\) for all \(a\in \overline{F_0}\) and \(\rho>0\), whereas \(\mu\) is locally finite on \(X\), as \(W\) is continuous; hence, \(\overline{F_0}\cap X=\emptyset\); as \(\overline{F_0}\subset\overline{X}\), we conclude that \(\overline{F_0}\subset\partial X\).

\medskip
\textit{Step~4: \(\overline{F_0}\) is a slice.}
We first collect some properties of \(d_0\).
As \(\overline{F_0}\cap X=\emptyset\) by Step~3, we have \(d_0>0\) on \(X\), whence \(f>0\) on \(X\), while \(f<1\) on \(X\), as \(s_{\max}=+\infty\); that is, \(X_{\mathrm{reg}}=X\).
Hence, by Step~1, \(d_0=\varphi-s_{\min}\) on the whole \(X\), so that, by~\eqref{eq:grad-phi},
\begin{equation}
\label{eq:nablad=N}
\nabla d_0=\nabla\varphi=\mathbf{N},
\end{equation}
which is smooth on \(X\); in particular, \(d_0\in C^2(X)\).
Moreover, for every \(y\in X\), by~\eqref{eq:y-pi0}, the point \(\pi_0(y)\in F_0\subset\overline{F_0}\) satisfies \(\|y-\pi_0(y)\|=d_0(y)=\dist(y;\overline{F_0})\), that is, it is a nearest point to \(y\) in \(\overline{F_0}\); as \(\overline{F_0}\) is closed and convex by Step~2, it is the unique one, that is, \(\pi_0(y)\) is the metric projection of \(y\) onto \(\overline{F_0}\).
In particular, \(\pi_0\) is the restriction to \(X\) of a \(1\)-Lipschitz map of \(\R^N\) which fixes \(\overline{F_0}\), and, by~\eqref{eq:y-pi0} and~\eqref{eq:nablad=N},
\begin{equation}
\label{eq:projection}
y-\pi_0(y)=d_0(y)\,\nabla d_0(y),
\qquad\forall\,y\in X.
\end{equation}

Let now \(L_0\) be the linear subspace parallel to \(L\), the affine hull of \(\overline{F_0}\), let \(Q\) be the orthogonal projection onto \(L_0^\perp\), let \(\ri(\overline{F_0})\) be the relative interior of \(\overline{F_0}\) in \(L\), which is nonempty, as \(\overline{F_0}\) is nonempty and convex, and let
\[
A=\{\,y\in X\,:\,\pi_0(y)\in\ri(\overline{F_0})\,\}.
\]
Since \(\pi_0\) is continuous and \(\ri(\overline{F_0})\) is relatively open in \(L\), the set \(A\) is open.
It is also nonempty: if \(b\in\ri(\overline{F_0})\subset\overline{X}\) and \(y_j\in X\) converge to \(b\), then \(\|\pi_0(y_j)-b\|\le\|y_j-b\|\to0\), as the metric projection onto \(\overline{F_0}\) is \(1\)-Lipschitz and fixes \(b\), so that \(y_j\in A\) for \(j\) large.

We first show that \(d_0=\dist(\,\cdot\,;L)\) on \(A\).
Let \(y\in A\). 
By~\eqref{eq:support-half-spaces}--\eqref{eq:projection}, we have that \(\langle z-\pi_0(y),y-\pi_0(y)\rangle\le0\) for all \(z\in \overline{F_0}\).
As \(\pi_0(y)\in\ri(\overline{F_0})\), for every \(v\in L_0\) we have \(\pi_0(y)\pm\varepsilon v\in \overline{F_0}\) for \(\varepsilon>0\) small enough; choosing \(z=\pi_0(y)\pm \varepsilon v\), we get \(y-\pi_0(y)\perp L_0\), so that \(\pi_0(y)\) is the orthogonal projection of \(y\) onto \(L\), and \(\dist(y;L)=\|y-\pi_0(y)\|=d_0(y)\).

Fix now \(p\in L\).
Every \(y\in\R^N\) can be written as \(y=p+v+Q(y-p)\), for some \(v\in L_0\), so that \(\dist(y;L)=\|Q(y-p)\|\).
Hence, as \(d_0=\dist(\,\cdot\,;L)\) on the open set \(A\), where \(d_0>0\), and as \(Q\) is symmetric and satisfies \(Q^T Q = Q^2 = Q\), we have
\begin{equation}
\label{eq:projection-2}
\nabla d_0(y)
=
\frac{Q(y-p)}{d_0(y)}
\in L_0^\perp,
\qquad
\forall\,y\in A.
\end{equation}
Differentiating once more and using the above equality, we get
\begin{equation}
\label{eq:hessian-cylinder}
D^2d_0(y)=\frac{Q-\nabla d_0(y)\otimes\nabla d_0(y)}{d_0(y)},
\qquad
\forall\,y\in A.
\end{equation}

We now show that \(A=X\).
Assume by contradiction that \(A\neq X\).
As \(X\) is connected and \(A\) is open and nonempty, there exists \(\xi\in X\setminus A\) in the closure of \(A\).
As \(d_0\in C^2(X)\) and \(d_0>0\) on \(X\), both sides of~\eqref{eq:projection-2} and~\eqref{eq:hessian-cylinder} are continuous on \(X\), so that these identities also hold at \(\xi\).
Consider now \(\pi_0(\xi)\in \overline{F_0}\setminus\ri(\overline{F_0})\), so that \(\xi=\pi_0(\xi)+d_0(\xi)\nabla d_0(\xi)\) by~\eqref{eq:projection}, and \(\nabla d_0(\xi)\perp L_0\) by~\eqref{eq:projection-2}.
As \(\pi_0(\xi)\) lies on the relative boundary of \(\overline{F_0}\) in \(L\), there exists \(\eta\in L_0\cap\mathbb{S}^{N-1}\) such that
\begin{equation}
\label{eq:support-K-in-L}
\langle z-\pi_0(\xi),\eta\rangle\le0
\qquad
\text{for all }z\in \overline{F_0}.
\end{equation}
For \(s\ge0\), let
\begin{equation}
\label{eq:y_s}
y_s=\pi_0(\xi)+d_0(\xi)\,\frac{\nabla d_0(\xi)+s\eta}{\sqrt{1+s^2}},
\end{equation}
which belongs to \(X\) for \(s\) small, as \(X\) is open and \(y_0=\xi\in X\).
Further, since \(z-\pi_0(\xi)\in L_0\) for all \(z\in \overline{F_0}\), by~\eqref{eq:projection-2} and~\eqref{eq:support-K-in-L}, we have
\[
\langle z-\pi_0(\xi),\nabla d_0(\xi)+s\eta\rangle\le0.
\]
By~\eqref{eq:y_s}, this is equivalent to \(\langle z-\pi_0(\xi),y_s-\pi_0(\xi)\rangle\le0\) for all \(z\in \overline{F_0}\), as \(y_s-\pi_0(\xi)\) is a positive multiple of \(\nabla d_0(\xi)+s\eta\).
Therefore, as
\begin{align*}
\|y_s-z\|^2
&
=\|y_s-\pi_0(\xi)\|^2+\|z-\pi_0(\xi)\|^2-2\langle z-\pi_0(\xi),y_s-\pi_0(\xi)\rangle
\\
&
\ge\|y_s-\pi_0(\xi)\|^2
\end{align*}
for all \(z\in \overline{F_0}\), we have \(\pi_0(y_s)=\pi_0(\xi)\).
Hence, as \(\nabla d_0(\xi)\perp\eta\) are unit vectors, by~\eqref{eq:y_s} we get \(d_0(y_s)=\|y_s-\pi_0(\xi)\|=d_0(\xi)\), so that, by~\eqref{eq:projection},
\[
y_s-\pi_0(\xi)
=
d_0(\xi)\,\nabla d_0(y_s).
\]
Differentiating at \(s=0\), as \(y_0=\xi\) and \(\frac{\mathrm{d}}{\mathrm{d}s}y_s\big|_{s=0}=d_0(\xi)\,\eta\) by~\eqref{eq:y_s}, we get
\[
d_0(\xi)\,\eta=d_0(\xi)^2\,D^2d_0(\xi)\,\eta.
\]
The left-hand side is nonzero, as \(d_0(\xi)>0\) and \(\eta\in\mathbb{S}^{N-1}\), whereas the right-hand side vanishes by~\eqref{eq:hessian-cylinder}, as \(Q\eta=0\) and \(\langle\nabla d_0(\xi),\eta\rangle=0\), since \(\eta\in L_0\).
This contradiction shows that \(A=X\), that is, \(d_0=\dist(\,\cdot\,;L)\) on \(X\), hence on \(\overline{X}\) by continuity.
It now immediately follows that \(L\cap\overline{X}\subset \overline{F_0}\), and thus \(\overline{F_0}=L\cap\overline{X}\), as \(\overline{F_0}\subset L\cap\overline{X}\).
By Step~3, \(L\cap X=\overline{F_0}\cap X=\emptyset\); in particular, \(L\neq\R^N\), that is, \(\dim L\le N-1\).

Finally, if \(L\) is a hyperplane, then \(X\), being connected and disjoint from \(L\), lies on one side of it, so that \(\mathbf{N}=\nabla d_0=\nabla\dist(\,\cdot\,;L)\) is constant on \(X\).
Conversely, if \(\mathbf{N}\equiv\theta\) on \(M_{c_0}\), by \cref{cor:parallel}\ref{item:par-half}, \(M_{c_0}\) is a relatively open subset of a hyperplane orthogonal to \(\theta\), so that by definition~\eqref{eq:F0}, \(F_0\) is a translate of \(M_{c_0}\), and thus spans a hyperplane.
Hence, \(L\) is a hyperplane if and only if \(\mathbf{N}\) is constant on \(M_{c_0}\), and \(\dim L\le N-2\) otherwise.
\end{proof}

\begin{remark}[Convexity and flatness]
\label{rem:flat-smooth}
In the proof of \cref{prop:focal-c}, the convexity of \(\overline{F_0}\) relies on the fact that the lines of steepest ascent of \(f\) are half-lines contained in \(X\), and on the convexity of \(X\), see also \cref{rem:nonconvex}.
The inclusion \(\overline{F_0}\subset\partial X\) relies on the non-integrability of the density near \(F_0\), given by \cref{prop:non-integrable}.
Finally, the flatness of \(\overline{F_0}\), in the sense that \(\overline{F_0}=L\cap\overline{X}\), relies on the continuity of the second fundamental form of the level sets of \(f\).
\end{remark}

\subsection{Proof of the \texorpdfstring{\cref{thm:main}}{Main Theorem}}

By \cref{lem:sweep}, \(f\) satisfies~\ref{hp:SW}, so that \cref{thm:dist-rep,cor:parallel} apply, and the level sets of \(f\) are parallel hypersurfaces.
If \(s_{\max}<+\infty\), then~(b) holds by \cref{prop:focal-point}.
If \(s_{\max}=+\infty\), then, by \cref{cor:parallel}, either \(\mathbf{N}\) is constant on \(M_{c_0}\), and~(a) holds, or \(s_{\min}>-\infty\), and the superlevel sets of \(f\) are the complements in \(X\) of the closed metric neighborhoods of \(F_0\).
In the latter case, by \cref{prop:focal-c}, the closure of \(F_0\), which has the same metric neighborhoods, is \(L\cap\overline{X}\), where \(L\) is its affine hull, \(\dim L\le N-2\), as \(\mathbf{N}\) is not constant, \(L\cap X=\emptyset\), and \(\dist(\,\cdot\,;F_0)=\dist(\,\cdot\,;L)\) on \(X\), so that~(c) holds.
The remaining claims follow from \cref{prop:width,prop:non-integrable}. \qed

\begin{remark}[Balls, half-spaces, and cylinders]
\label{rem:focal}
The \cref{thm:main} is consistent with the examples of \cite[Sect.~6]{BS26}: concentric balls appear in \cite[Sect.~6.1]{BS26} for measures with \(\mu(X)=+\infty\), the measures of \cite[Sect.~6.2]{BS26} and the Gauss-like product measures of~\cite{BCM16}, for which parallel half-spaces occur, are finite, and in \cite[Sect.~6.4]{BS26} balls occur for finite measures, but they are isoperimetric only for small volumes, see \cref{rem:small-volumes}.
Also the third alternative does occur, with affine subspaces \(L\) of every dimension \(\ell\in\{0,\dots,N-2\}\), see \cref{sec:examples-c}, while parallel half-spaces also occur for infinite measures, as complements of tubes around a hyperplane, see \cref{rem:hyperplane}.

By the \cref{thm:main}, cylinders around affine subspaces of positive dimension never occur.
This is consistent with the fact that cylinders may well be isoperimetric for large volumes, see, \emph{e.g.},~\cite{FMP11,RV15b}, while the proof of \cref{prop:focal-point} only uses the isoperimetric property for small volumes.
On the contrary, complements of cylinders may occur, see \cref{ex:wedge}.
\end{remark}

\begin{remark}[Small volumes]
\label{rem:small-volumes}
By~\cite[Rem.~5.5]{BS26}, Talenti's inequality still holds for functions whose support has measure at most \(\overline{m}\in(0,\mu(X))\), if~\ref{hp:T1}--\ref{hp:T3} are replaced by their counterparts (T1\('\))--(T3\('\)) of~\cite[Sect.~2]{BS26}, which only involve volumes not exceeding \(\overline{m}\), as for the convex polytopes of~\cite[Sect.~6.4.1]{BS26}.
The proofs of \cref{lem:sweep,prop:focal-point} then apply to the superlevel sets of \(f\) of measure smaller than \(\overline{m}\), which are therefore mutually parallel, and concentric balls if \(s_{\max}<+\infty\); in particular, tubes around sets of positive dimension are still excluded.
Below this threshold, however, the lines of steepest ascent of \(f\) may leave \(X\), as in~\cite[Sect.~6.4]{BS26}, and we do not pursue this case.
\end{remark}

\section{Examples of alternative~(c)}
\label{sec:examples-c}

In this section, we show that alternative~(c) of the \cref{thm:main} does occur, with an affine subspace \(L\) of any dimension \(\ell\in\{0,\dots,N-2\}\), and that alternative~(a) also occurs for infinite measures, see \cref{rem:hyperplane}.
Let \(N\ge2\), let \(\ell\in\{0,\dots,N-2\}\), and let \(n=N-\ell\ge2\).
We write points of \(\R^N=\R^n\times\R^\ell\) as \(x=(y,z)\), with the convention that the variable \(z\) is absent if \(\ell=0\).
Let \(C\subset\R^n\), \(C\neq\R^n\), be an open convex cone with vertex at the origin, so that \(0\notin C\), let \(\Sigma=C\cap\mathbb{S}^{n-1}\), which is relatively open in \(\mathbb{S}^{n-1}\), and let \(\sigma_C=\mathcal{H}^{n-1}(\Sigma)\in(0,+\infty)\).
For \(k>0\), we let
\[
X=C\times(0,1)^\ell\subset\R^N,
\qquad
e^{W(y,z)}=\|y\|^{-n-k},
\qquad (y,z)\in X.
\]
Then \(X\) is open and convex, and \(e^W\) is of class \(C^1\) on \(X\), but not integrable near the face \(\{0\}\times[0,1]^\ell\subset\partial X\), which is a vertex if \(\ell=0\), and an edge if \(\ell=1\).
Given any \(r>0\), we let
\[
X_r=\{\,(y,z)\in X\,:\,\|y\|>r\,\},
\]
that is, the complement in \(X\) of the closed tube \(\overline{B_r(0)}\times\R^\ell\), where \(B_r(0)\) denotes a ball of \(\R^n\).
Integrating in polar coordinates in the variable \(y\), we have
\begin{equation}
\label{eq:complements-vol-per}
\mu(X_r)=\frac{\sigma_C}{k\, r^{k}},
\qquad
\Per_\mu(X_r;X)=\frac{\sigma_C}{r^{k+1}},
\end{equation}
while \(\mu(X\setminus X_r)=+\infty\).
In particular, for every \(m\in(0,+\infty)\), there exists exactly one \(r>0\) such that \(\mu(X_r)=m\).
In the next proposition, via a calibration argument, we show that the sets \(X_r\) are isoperimetric in this setting.

\begin{proposition}
\label{prop:iso-complements}
For every Borel set \(E\subset X\) with \(0<\mu(E)<+\infty\),
\[
\Per_\mu(E;X)
\ge
\Per_\mu(X_r; X),
\]
where \(r\) is such that \(\mu(E)=\mu(X_r)\).
\end{proposition}

\begin{proof}
Let \(E\subset X\) be a Borel set with \(0<\mu(E)<+\infty\) and \(\Per_\mu(E;X)<+\infty\), the claim being trivial otherwise.
By De Giorgi's structure theorem, see~\cite[Thm.~15.9]{Mag12book},
\begin{equation}
\label{eq:gauss-green}
\int_E\operatorname{div}\Phi\,\mathrm{d}x=\int_{\partial^*E\cap X}\langle\Phi,\nu_E\rangle\,\mathrm{d}\mathcal{H}^{N-1},
\qquad\forall\,\Phi\in C^1_c(X;\R^N),
\end{equation}
where \(\nu_E\) is the measure-theoretic outer unit normal to \(E\).
We let \(r>0\) be such that \(\mu(X_r)=\mu(E)\).

\medskip
\textit{Step~1: The calibration.}
Let \(h\) be given by
\[
h(s)
=
\begin{cases}
\dfrac{(k+1)s-s^{k+1}}{k}
\qquad
&s\in[0,1],
\\[5pt]
1,
&s>1.
\end{cases}
\]
This function is easily seen to be in \(C^1([0,+\infty))\), to be nondecreasing, and to satisfy
\begin{equation}
\label{eq:bound-h}
0\le h(s)\le\min\left\{1,\frac{k+1}{k}\, s\right\}.
\end{equation}
We now consider the vector field \(T\) defined on \(X\) as
\[
T(y,z)
=
-h\left(\frac{\|y\|}{r}\right)
\,e^{W(y,z)}
\,\left(\frac{y}{\|y\|},0\right),
\qquad (y,z)\in X.
\]
Clearly, \(\|T\|\le e^{W}\). %
Further, by definition of \(h\), recalling that \(e^W = \|y\|^{-n-k}\), and recalling that \(y\neq0\) for all \((y,z)\in X\), the vector field \(T\) is of class \(C^1\) on \(X\).
A straightforward computation, using that the divergence in \(\R^n\) of \(y/\|y\|\) equals \((n-1)/\|y\|\), shows that
\[
\Div T
=
\frac{1}{\|y\|^{n+k}\,r}\left(\frac{(k+1)h(s)}{s}-h'(s)\right),
\qquad
\text{where }
s = \frac{\|y\|}{r}.
\]
In particular, the term in brackets equals \(k+1\) if \(s\le1\), and \((k+1)/s\) if \(s>1\).
Thence, recalling once again that \(e^W = \|y\|^{-n-k}\), we have
\begin{equation}
\label{eq:div-T}
\Div T=\lambda\,e^W,
\qquad
\text{where }
\lambda(y,z)=\frac{k+1}{\max\{\|y\|,r\}}.
\end{equation}

\medskip
\textit{Step~2: Localization.}
We now wish to build from \(T\) a vector field \(\Phi\) admissible in~\eqref{eq:gauss-green}.
Therefore, we need to cut-off nearby the boundary of \(X\) and also localize in a bounded set.
Thus, we consider the following two cut-off functions: given \(0<\varepsilon<r<R\) and \(\delta>0\), we let
\[
\psi\in C^1_c((\varepsilon,2R))
\qquad
\text{and}
\quad
\zeta\in C^1_c(\Sigma\times(0,1)^\ell)
\]
taking both values in \([0,1]\), with
\begin{equation}
\label{eq:cut-offs}
\begin{gathered}
\psi=1
\text{ on }
[2\varepsilon,R],
\quad
0\le\psi'\le\frac2\varepsilon
\text{ on }
(\varepsilon,2\varepsilon),
\quad
\psi'\le0
\text{ on }
(R,2R),
\\
\zeta=1
\text{ on }
\Sigma_\delta\times[\delta,1-\delta]^\ell,
\,\,
\text{with }
\Sigma_\delta=\{\,\omega\in\Sigma\,:\,\dist(\omega;\mathbb{S}^{n-1}\setminus\Sigma)\ge\delta\,\}.
\end{gathered}
\end{equation}
We let \(\psi\) act on the \(\|y\|\) variable, thus providing a localization in the set \((B_{2R}(0)\setminus \overline{B_{\varepsilon}(0)})\times(0,1)^\ell\), and we let \(\zeta\) act on the angular variable \(y/\|y\|\) and on \(z\), thus killing \(T\) near the lateral walls \(\partial C\times(0,1)^\ell\) and near the faces \(C\times\partial(0,1)^\ell\) of \(X\).
Namely, we let
\[
\Phi(y,z)=\psi(\|y\|)\,\zeta\left(\frac{y}{\|y\|},z\right)\,T(y,z).
\]
By definition it belongs to \(C^1_c(X;\R^N)\), with \(\|\Phi\|\le e^W\).

\medskip
\textit{Step~3: Comparison with \(X_r\).}
Since \(y\mapsto y/\|y\|\) is constant along rays from the origin, the gradient in \(y\) of \(\zeta(y/\|y\|,z)\) is orthogonal to \(y\); as \(T\) has no \(z\)-component and its \(y\)-component is parallel to \(y\), the gradient of \(\zeta(y/\|y\|,z)\) is orthogonal to \(T\), so that, by~\eqref{eq:div-T} and the definition of \(T\),
\begin{align*}
\Div\Phi
&
=
\psi\,\zeta\,\Div T
+
\zeta\,\bigl\langle\nabla(\psi(\|y\|)),T\bigr\rangle
+
\psi\,\Bigl\langle\nabla\Bigl(\zeta\Bigl(\frac{y}{\|y\|},z\Bigr)\Bigr),T\Bigr\rangle
\\
&
=
\zeta\,\bigl(\psi\,\Div T+\bigl\langle\nabla(\psi(\|y\|)),T\bigr\rangle\bigr)
\\
&
=
\zeta(\psi \lambda  -  \psi' h ) e^W,
\end{align*}
where \(\psi,\psi'\) are evaluated at \(\|y\|\), \(\zeta\) at \((y/\|y\|,z)\), and \(h\) at \(\|y\|/r\).

Since \(\|\Phi\|\le e^W\), by the previous equality, by~\eqref{eq:gauss-green}, and recalling that \(\mu = e^W \mathcal{L}^N\), we get
\[
\int_E \zeta \psi \lambda \,\mathrm{d}\mu
\le
\Per_\mu (E; X)
+
\int_E \zeta \psi' h \,\mathrm{d}\mu.
\]
Owing to~\eqref{eq:cut-offs}, \(\psi'\) is supported on \([\varepsilon, 2\varepsilon]\cup[R, 2R]\) and in the interval \([R, 2R]\) we have \(\psi' h\le 0\); recalling also that \(\zeta\) is nonnegative, we can bound the second term on the right-hand side by considering only the integral over the set \(E\cap\{\varepsilon <\|y\|<2\varepsilon\}\); using that \(\zeta\le 1\), and the upper bound~\eqref{eq:bound-h} on \(h\) and that on \(\psi'\) in~\eqref{eq:cut-offs}, we get
\begin{align*}
\int_E \zeta \psi' h\,\mathrm{d}\mu
&
\le
\int_{E\cap \{\varepsilon < \|y\| <2\varepsilon\}} \psi'(\|y\|) h\left( \frac{\|y\|}{r}\right)\,\mathrm{d}\mu
\\
&
\le
\int_{E\cap \{\varepsilon < \|y\| <2\varepsilon\}} \frac2\varepsilon \frac{k+1}{rk} \|y\|\,\mathrm{d}\mu
\\
&
\le
\frac{4(k+1)}{kr}
\mu(E\cap \{\varepsilon < \|y\| <2\varepsilon\}).
\end{align*}
Thence, using this latter estimate in the previous one, we have the upper bound
\[
\int_E \zeta \psi \lambda \,\mathrm{d}\mu
\le
\Per_\mu (E; X)
+
\frac{4(k+1)}{kr}
\mu(E\cap \{ \|y\| <2\varepsilon\}).
\]
Letting \(\delta,\varepsilon\to0^+\) and \(R\to+\infty\), owing to the dominated convergence theorem, as \(\lambda\le(k+1)/r\) by~\eqref{eq:div-T} and \(\mu(E)<+\infty\), we conclude that
\[
\Per_\mu(E;X)
\ge
\int_E\lambda\,\mathrm{d}\mu
=
\int_{E\cap X_r}\lambda\,\mathrm{d}\mu
+
\int_{E\setminus X_r}\lambda\,\mathrm{d}\mu.
\]

By definition of \(\lambda\), see~\eqref{eq:div-T}, and of \(X_r\), we have that
\begin{equation*}
\lambda=\frac{k+1}{r}
\quad
\text{ on } E\setminus X_r,
\qquad\qquad
\lambda< \frac{k+1}{r}
\quad \text{ on } X_r\setminus E;
\end{equation*}
and since \(\mu(E\setminus X_r)=\mu(X_r\setminus E)<+\infty\), we can bound from below the second term with the integral of \(\lambda\) over \(X_r\setminus E\), so that, using the definition of \(\lambda\) and of \(X_r\), that of \(e^W\), passing to polar coordinates in the variable \(y\), and using~\eqref{eq:complements-vol-per}, we get
\[
\Per_\mu(E; X)
\ge
\int_{X_r} \lambda\,\mathrm{d}\mu
=
(k+1)\,\sigma_C \int_{r}^{+\infty} \frac{1}{\rho^{k+2}}\,\mathrm{d}\rho
=
\Per_\mu(X_r; X).
\qedhere
\]
\end{proof}

\begin{example}[Complements of tubes around a face]
\label{ex:wedge}
Let \(X\) and \(\mu\) be as above, and let
\[
f(y,z)=\frac{\|y\|^2}{1+\|y\|^2},
\qquad (y,z)\in X.
\]
Then,~\ref{hp:N1} and~\ref{hp:N2} hold, as \(\inf f=0\) and \(\sup f=1\) are not attained in \(X\), and \(\nabla f\neq0\) in \(X\).
The superlevel sets of \(f\) are the sets \(X_r\), which have finite measure by~\eqref{eq:complements-vol-per}, so that~\ref{hp:N3} holds, and which saturate the isoperimetric inequality of \cref{prop:iso-complements}, so that~\ref{hp:T1} and~\ref{hp:T2} hold with \(I(0)=0\) and \(I(m)=\Per_\mu(X_r;X)\) for \(\mu(X_r)=m>0\), that is, by~\eqref{eq:complements-vol-per},
\[
I(m)=k^{\frac{k+1}{k}}\sigma_C^{-\frac1k}\,m^{\frac{k+1}{k}},
\]
which is continuous.
Its level sets are the sets
\[
(\partial B_r(0)\cap C)\times(0,1)^\ell=r\Sigma\times(0,1)^\ell,
\]
which are connected.
Its lines of steepest ascent are the half-lines emanating from the face \(\{0\}\times(0,1)^\ell\) given by \(\{\,(\tau\omega,z)\,:\,\tau>0\,\}\), with \(\omega\in\Sigma\) and \(z\in(0,1)^\ell\).
Along these, \(f\) takes all values in \((0,1)\).
Hence, all the assumptions of the \cref{thm:main} hold, \(\mathbf{N}\) is not constant, \(s_{\min}>-\infty\), and \(s_{\max}=+\infty\), so that we are in alternative~(c).
In particular, the closure of the focal set is \(\{0\}\times[0,1]^\ell=L\cap\overline{X}\), with \(L=\{0\}\times\R^\ell\) and \(\dim L=\ell\), in accordance with \cref{prop:focal-c}.
Further, as prescribed by \cref{prop:non-integrable}, the density is not integrable near any point of \(F_0\).
Finally, we notice that the P\'olya--Szeg\H{o} principle of~\cite{BS26} applies to \(f\), and so do the Faber--Krahn, Saint-Venant, and \((p,q)\)-eigenvalue inequalities that it implies, see~\cite[Sect.~6]{BS26}, with the sets \(X_r\) as extremal sets.
\end{example}

\begin{remark}
For \(\ell=0\) and \(N=2\), the set \(X=C\) is a planar convex sector, and we recover that the complements of the disks centered at the vertex are isoperimetric in \(C\) endowed with the density \(\|y\|^{-2-k}\), see also~\cite[Remark after Thm.~4.17]{DHHT12} and~\cite[Prop.~4.3]{CJQW08}; this is the example of~\cite[Sect.~6.3]{BS26}, for which alternative~(c) holds with \(F_0=\{0\}\).
For \(\ell=0\) and \(N\ge3\), see also the remark after~\cite[Prop.~7.5]{DHHT12}.
\end{remark}

\begin{remark}[Half-spaces as complements of tubes]
\label{rem:hyperplane}
The above construction also makes sense for \(\ell=N-1\), that is, \(n=1\), with \(C=(0,+\infty)\), so that \(\Sigma=\{1\}\), \(\sigma_C=1\), and \(X=(0,+\infty)\times(0,1)^{N-1}\) is a half-cylinder whose cross-section is the hypercube \((0,1)^{N-1}\).
In this case, the divergence in \(\R\) of \(y/\|y\|\) vanishes, consistently with the formula \((n-1)/\|y\|\), and the proof of \cref{prop:iso-complements} applies verbatim, showing that the sets
\[
X_r=\{\,(y,z)\in X\,:\,y>r\,\},
\]
that is, the intersections of \(X\) with the half-spaces orthogonal to its axis, are isoperimetric.
As in~\cite{BCM16}, half-spaces are thus isoperimetric in a cylinder endowed with a product density; there, however, \(\mu\) is finite and Gauss-like in the axial direction, whereas here \(\mu(X)=+\infty\) and the density \(y^{-1-k}\) blows up at the face \(\{0\}\times[0,1]^{N-1}\), that is, at the base of the half-cylinder.

Therefore, as in \cref{ex:wedge}, the function \(f(y,z)=y^2/(1+y^2)\) satisfies all the assumptions of the \cref{thm:main}, with \(s_{\min}>-\infty\) and \(s_{\max}=+\infty\), but now \(\mathbf{N}\) is constant, and we are in alternative~(a): the closure of the focal set is the face \(\{0\}\times[0,1]^{N-1}\), whose affine hull is the hyperplane \(\{0\}\times\R^{N-1}\).
This shows that, in \cref{prop:focal-c}, the bound \(\dim L\le N-2\) does require \(\mathbf{N}\) not to be constant.
\end{remark}

The calibration of \cref{prop:iso-complements} also shows that the convexity of \(X\) cannot be dropped in the \cref{thm:main}, even though Talenti's theorem only requires \(X\) to be open and connected.

\begin{remark}[Convexity cannot be dropped]
\label{rem:nonconvex}
Let \(K\subset\R^N\) be a convex body with smooth boundary, and let \(X=\R^N\setminus K\), which is open and connected, but not convex.
Every \(x\in X\) can be uniquely written as \(x=p+s\,\nu_K(p)\), with \(p\in\partial K\), \(s=d(x)=\dist(x;K)\), and \(\nu_K\) the outer unit normal to \(K\); in these coordinates, the Lebesgue measure reads
\[
J\,\mathrm{d}s\,\mathrm{d}\mathcal{H}^{N-1}(p),
\quad
\text{with }
J(p,s)=\prod_{i=1}^{N-1}(1+s\,\kappa_i(p)),
\]
where \(\kappa_i\ge0\) are the principal curvatures of \(\partial K\), and \(\Delta d=\partial_s\log J\).
For \(k>0\), let
\[
e^W=\frac{d^{-1-k}}{J},
\qquad
f=\frac{d^2}{1+d^2},
\]
so that \(\mu\) is \(s^{-1-k}\,\mathrm{d}s\,\mathrm{d}\mathcal{H}^{N-1}(p)\), and the superlevel sets of \(f\) are the sets \(X_r\) defined as \(\{d>r\}\).
With \(h\) as in the proof of \cref{prop:iso-complements}, the field \(T=-h(d/r)\,e^W\,\nabla d\) is of class \(C^1\) on \(X\).
As \(\partial K\) is compact, the proof of \cref{prop:iso-complements} applies verbatim, with \(\|y\|\) replaced by \(d\) and with no cut-off \(\zeta\), and shows that the sets \(X_r\) are isoperimetric.
Hence, as in \cref{ex:wedge}, \(f\) satisfies~\ref{hp:N1}--\ref{hp:N3} and~\ref{hp:T1}--\ref{hp:T4}.

However, the sets \(X_r\) are complements of the parallel bodies \(K\oplus\overline{B_r(0)}\) and fall into none of the alternatives of the \cref{thm:main}: if \(K\) is not a ball, their boundaries are neither hyperplanes, nor spheres, nor cylinders, and if \(K\) is a ball, they are complements of balls whose center does not lie in \(\overline{X}\).
In the proof of the \cref{thm:main}, what breaks down is Step~2 of the proof of \cref{prop:focal-c}: \(f\) is a function of the distance from \(\partial K\), which plays the role of \(F_0\), but the convex hull of \(\partial K\), that is, \(K\), is not contained in \(\overline{X}\).
\end{remark}

\bibliographystyle{plainurl}
\bibliography{talenti-geometries}

\begin{thebibliography}{10}

\bibitem{AL26}
P.~Acampora and J.~Lamboley.
\newblock Sharp quantitative {T}alenti inequality in particular cases.
\newblock {\em Arch. Ration. Mech. Anal.}, 250:Paper No. 24, 2026.
\newblock \href {https://doi.org/10.1007/s00205-026-02179-3}
  {\path{doi:10.1007/s00205-026-02179-3}}.

\bibitem{BS26}
G.~Bartoli and G.~Saracco.
\newblock Spectral inequalities for weighted {$p$}-{L}aplacians via {T}alenti
  symmetrization, 2026.
\newblock \href {https://arxiv.org/abs/2605.29721} {\path{arXiv:2605.29721}}.

\bibitem{Bol73}
J.~Bolton.
\newblock Transnormal systems.
\newblock {\em Quart. J. Math. Oxford Ser. (2)}, 24(1):385--395, 1973.
\newblock \href {https://doi.org/10.1093/qmath/24.1.385}
  {\path{doi:10.1093/qmath/24.1.385}}.

\bibitem{BCM12}
F.~Brock, F.~Chiacchio, and A.~Mercaldo.
\newblock Weighted isoperimetric inequalities in cones and applications.
\newblock {\em Nonlinear Anal.}, 75(15):5737--5755, 2012.
\newblock \href {https://doi.org/10.1016/j.na.2012.05.011}
  {\path{doi:10.1016/j.na.2012.05.011}}.

\bibitem{BCM14}
F.~Brock, F.~Chiacchio, and A.~Mercaldo.
\newblock A weighted isoperimetric inequality in an orthant.
\newblock {\em Potential Anal.}, 41(1):171--186, 2014.
\newblock \href {https://doi.org/10.1007/s11118-013-9367-4}
  {\path{doi:10.1007/s11118-013-9367-4}}.

\bibitem{BCM16}
F.~Brock, F.~Chiacchio, and A.~Mercaldo.
\newblock An isoperimetric inequality for {G}auss-like product measures.
\newblock {\em J. Math. Pures Appl. (9)}, 106(2):375--391, 2016.
\newblock \href {https://doi.org/10.1016/j.matpur.2016.02.014}
  {\path{doi:10.1016/j.matpur.2016.02.014}}.

\bibitem{BZ88}
J.~E. Brothers and W.~P. Ziemer.
\newblock Minimal rearrangements of {S}obolev functions.
\newblock {\em J. Reine Angew. Math.}, 384:153--179, 1988.
\newblock \href {https://doi.org/10.1515/crll.1988.384.153}
  {\path{doi:10.1515/crll.1988.384.153}}.

\bibitem{CROS16}
X.~Cabr\'{e}, X.~Ros-Oton, and J.~Serra.
\newblock Sharp isoperimetric inequalities via the {ABP} method.
\newblock {\em J. Eur. Math. Soc. (JEMS)}, 18(12):2971--2998, 2016.
\newblock \href {https://doi.org/10.4171/JEMS/659}
  {\path{doi:10.4171/JEMS/659}}.

\bibitem{CJQW08}
C.~Carroll, A.~Jacob, C.~Quinn, and R.~Walters.
\newblock The isoperimetric problem on planes with density.
\newblock {\em Bull. Austral. Math. Soc.}, 78(2):177--197, 2008.
\newblock \href {https://doi.org/10.1017/S000497270800052X}
  {\path{doi:10.1017/S000497270800052X}}.

\bibitem{DHHT12}
A.~D\'{\i}az, N.~Harman, S.~Howe, and D.~Thompson.
\newblock Isoperimetric problems in sectors with density.
\newblock {\em Adv. Geom.}, 12(4):589--619, 2012.
\newblock \href {https://doi.org/10.1515/advgeom-2012-0009}
  {\path{doi:10.1515/advgeom-2012-0009}}.

\bibitem{Ehr83}
A.~Ehrhard.
\newblock Sym\'etrisation dans l'espace de {G}auss.
\newblock {\em Math. Scand.}, 53(2):281--301, 1983.
\newblock \href {https://doi.org/10.7146/math.scand.a-12035}
  {\path{doi:10.7146/math.scand.a-12035}}.

\bibitem{FMP11}
N.~Fusco, F.~Maggi, and A.~Pratelli.
\newblock On the isoperimetric problem with respect to a mixed
  {E}uclidean--{G}aussian density.
\newblock {\em J. Funct. Anal.}, 260(12):3678--3717, 2011.
\newblock \href {https://doi.org/10.1016/j.jfa.2011.01.007}
  {\path{doi:10.1016/j.jfa.2011.01.007}}.

\bibitem{Mag12book}
F.~Maggi.
\newblock {\em {S}ets of {F}inite {P}erimeter and {G}eometric {V}ariational
  {P}roblems}, volume 135 of {\em Cambridge Studies in Advanced Mathematics}.
\newblock Cambridge University Press, Cambridge, 2012.
\newblock An {I}ntroduction to {G}eometric {M}easure {T}heory.
\newblock \href {https://doi.org/10.1017/CBO9781139108133}
  {\path{doi:10.1017/CBO9781139108133}}.

\bibitem{Mag91}
R.~Magnanini.
\newblock On functions whose curves of steepest descent are straight lines.
\newblock {\em Appl. Anal.}, 41(1-4):171--182, 1991.
\newblock \href {https://doi.org/10.1080/00036819108840022}
  {\path{doi:10.1080/00036819108840022}}.

\bibitem{PS51book}
G.~P\'{o}lya and G.~Szeg{\H{o}}.
\newblock {\em Isoperimetric {I}nequalities in {M}athematical {P}hysics}.
\newblock Annals of Mathematics Studies, No. 27. Princeton University Press,
  Princeton, NJ, 1951.

\bibitem{RV15b}
M.~Ritor\'{e} and E.~Vernadakis.
\newblock Isoperimetric inequalities in convex cylinders and cylindrically
  bounded convex bodies.
\newblock {\em Calc. Var. Partial Differential Equations}, 54(1):643--663,
  2015.
\newblock \href {https://doi.org/10.1007/s00526-014-0800-3}
  {\path{doi:10.1007/s00526-014-0800-3}}.

\bibitem{Sch90book}
H.~A. Schwarz.
\newblock {\em Gesammelte mathematische {A}bhandlungen. {B}and {I}, {II}}.
\newblock Chelsea Publishing Co., Bronx, NY, 1972.
\newblock Nachdruck in einem Band der Auflage von 1890.
\newblock \href {https://doi.org/10.1007/978-3-642-50665-9}
  {\path{doi:10.1007/978-3-642-50665-9}}.

\bibitem{Seg38}
B.~Segre.
\newblock Famiglie di ipersuperficie isoparametriche negli spazi {E}uclidei ad
  un qualunque numero di dimensioni.
\newblock {\em Atti Accad. Naz. Lincei Rend. Cl. Sci. Fis. Mat. Nat. (6)},
  27:203--207, 1938.

\bibitem{Tal97}
G.~Talenti.
\newblock A weighted version of a rearrangement inequality.
\newblock {\em Ann. Univ. Ferrara Sez. VII}, 43:121--133, 1997.
\newblock \href {https://doi.org/10.1007/BF02837230}
  {\path{doi:10.1007/BF02837230}}.

\bibitem{Wan87}
Q.-M. Wang.
\newblock Isoparametric functions on {R}iemannian manifolds. {I}.
\newblock {\em Math. Ann.}, 277(4):639--646, 1987.
\newblock \href {https://doi.org/10.1007/BF01457863}
  {\path{doi:10.1007/BF01457863}}.

\end{thebibliography}

\end{document}